\documentclass[12pt]{article}
\usepackage{amsmath,subfigure}
\usepackage{amssymb}
\usepackage{epsfig}
\usepackage{eepic}
\usepackage{graphicx}
\usepackage{float}
\usepackage{latexsym,graphics,color}
\usepackage{dsfont}
\newtheorem{theorem}{Theorem}[section]
\newtheorem{proposition}[theorem]{Proposition}
\newtheorem{lemma}[theorem]{Lemma}

\newtheorem{definition}[theorem]{Definition}

\newtheorem{assumption}[theorem]{Assumption}

\newtheorem{remark}{Remark}[section]

\def\sqr#1#2{{\vcenter{\vbox{\hrule height.#2pt\hbox{\vrule width.#2pt height#1pt \kern#1pt\vrule width.#2pt}\hrule height.#2pt}}}}

\newcommand{\cR}{\mathcal{R}}

\begin{document}
  \vspace{1.5cm}
                   \begin{center}
   {\large \bf Distributed space scales in a semilinear reaction-diffusion system including a parabolic variational inequality:  A well-posedness  study }\\[1cm]
                   \end{center}
\begin{center}
Dedicated to Professor Marek Niezg\'odka on his 60th birthday
\end{center}

                   \begin{center}
                  {\sc Tasnim Fatima$^{\rm \dagger}$, Adrian Muntean$^{\rm \dagger,\star}$ and Toyohiko Aiki$^{\rm \ddag}$}\\
                    \vspace{0.3cm}
$^{\rm \dagger,\star}$CASA - Centre for Analysis,
Scientific computing and Applications, Department of Mathematics and
Computer Science, Eindhoven University of Technology, P.O. Box 513,
5600 MB Eindhoven, The Netherlands.  e-mail: t.fatima@tue.nl
                   \vspace{0.3cm}

$^{\rm \star}$Institute for Complex Molecular Systems (ICMS),
 e-mail: a.muntean@tue.nl
  \vspace{0.3cm}

 $^{\rm \ddag}$Department of Mathematics,
Faculty of Education,
Gifu University
Yanagido 1-1, Gifu, 501-1193
Japan. e-mail: aiki@gifu-u.ac.jp
%
       \end{center}
       \vspace{1cm}

\hspace*{-0.6cm}{\bf Abstract.} This paper treats the solvability of a semilinear reaction-diffusion system,
 which incorporates transport (diffusion) and reaction effects emerging from two separated spatial scales:
  $x$ - macro and $y$ - micro. The system's origin connects to the  modeling of concrete corrosion in sewer
  concrete pipes. It consists of three partial differential equations which are mass-balances
   of concentrations, as well as, one ordinary differential equation tracking the damage-by-corrosion.
The system is semilinear, partially dissipative, and coupled
via the solid-water interface at the microstructure (pore) level.
The  structure of the model equations is obtained in \cite{tasnim1} by upscaling of
the physical and chemical processes taking place within the microstructure of the concrete.
Herein we ensure the positivity and  $L^\infty-$bounds on concentrations, and then prove
the global-in-time existence and uniqueness of a suitable class of positive and bounded
solutions that are stable with respect to the two-scale data and model parameters.
The main ingredient to prove existence include fixed-point arguments and
convergent two-scale Galerkin approximations.  \\

\\ 

 \section{Introduction}

We consider a two-scale (distributed-microstructure\footnote{This terminology is very much due to R. E Showalter; see  chapter 9 in \cite{hornung}.}- or double-porosity-) system modeling penetration of corrosion in concrete sewer pipes.
This kind of models appears in a multitude of real-world applications and are therefore of great importance
 mainly because they are able to connect the information from the microscale to the  macroscale (e.g. via the boundary of the cell, micro-macro transmission conditions).
 They are usually obtained in the homogenization limit as the scale of the inhomogeneity goes to zero.
These models provide a way to represent a continuous distribution of cells within a global reference geometry.
Roughly speaking, to each point $x\in\Omega$, we assign a representative cell $Y_x$. The flow within
 each cell is described (independently w.r.t what happens at the macroscale) by an initial-boundary-value problem.
The solution of the problem posed in the cell $Y_x$ is coupled via the boundary of $Y_x$ to the macroscale.
In order to derive such models, different choice of microstructures within the global
domain can be considered. Note that the microstructure $Y_x$ does not necessarily need to be periodically distributed in $\Omega$.
Uniformly periodic (i.e., $Y_x=Y$) and locally periodic are some of the options. Apparently, this sort of models are good approximation of
real situations, but involve the computation of a very large number of cells problems.

The model was obtained via formal homogenization using the different scalings in $\epsilon$ (a scale parameter referring to the micro-geometry) of the diffusion coefficients in locally periodic case;
details can be seen in \cite{tasnim1}.
The relevant mathematical questions at this point are threefold:
\begin{itemize}
 \item[(Q1)] How does the information 'flow' between macro and micro scales? What are the correct micro-macro transmission boundary conditions? Is the model  well-posed in the sense of Hadamard?
 \item[(Q2)] What is the long-time behavior of this two-scale reaction-diffusion system?
 \item[(Q3)] To which extent is the model computable?
 \end{itemize}

In this paper, we focus on the first type of questions -- the well-posedness of the two-scale model --
 preparing in the same time the playground to tackle the second question.
For more information on the modeling, analysis and simulation of two-scale scenarios, we refer the reader to
\cite{tasnim2,Meier,Meier1,maria1,maria2,Treutler,Noorden1}.
 We postpone the study of the long-time behavior (and of the inherently occurring memory effects) to a later stage.
  Also, we will investigate elsewhere the computability of our system. However, it is worth mentioning
   now preliminary results in this direction: \cite{Lakkis} proves the rate of convergence
    for a two-scale Galerkin scheme, while the authors of \cite{tasnim2} we produce the  first numerical simulations of two-scale effects to our system.
    Our working techniques are inspired by \cite{maria1,maria2}. The novelty we bring in this paper is two-fold:
    \begin{itemize}
      \item We are able to apply the two-scale Galerkin procedure to a partly-dissipative R-D system, (compare to the case of fully-dissipative system treated in \cite{maria1,maria2}).
      \item We are able to circumvent the use of both extra (macro) $x-$regularity of the micro-solutions and $H^2-$regularity w.r.t. $y$ of the initial data by making use of the parabolic variational inequality framework for part of our system.
    \end{itemize}
The main reason why we want to keep low the regularity assumption on the initial data is that, at a later stage, we would like to couple this reaction-diffusion system with the actual degradation of the material, namely with partial differential  equations governing the mechanics of fracturing concrete.

The paper is organized as follows: Section \ref{pro} contains a brief introduction to the chemistry of the problem as well as a concise  description of the geometry.
The two-scale model equations are introduced  in Section \ref{Model} while Section \ref{setting} includes
 the functional setting  and assumptions as well as the  main results of the paper.
 The proofs of the results are given in Section \ref{Proof}.

\section{Description of the problem}\label{pro}

\subsection{Chemistry}
Sulfuric attack to concrete structures is one of the most aggressive chemical attacks.
This happens usually in sewer pipes, since in sewerage there is a lot of
hydrogen sulfide. Anaerobic bacteria present in the waste flow produce hydrogen sulfide gaseous $H_2S$.
 From the air space of the pipe, $H_2S$ enters the air space of the microstructure where it diffuses and dissolves in the pore water. Here aerobic bacteria
catalyze into sulfuric acid $H_2SO_4$. As a next step, $H_2SO_4$ reacts with calcium carbonate
(part of solid matrix, say concrete) and eventually destroys  locally the pipe by spalling. The model we study here incorporates two particular mechanisms:
\begin{itemize}
  \item the exchange of $H_2S$ in water and air phase of the microstructure and vise versa \cite{Balls},
  \item production of $gypsum$ as a result of reaction between $H_2SO_4$ and calcium carbonate at solid-water interface.
\end{itemize}
   The transfer of $H_2S$ between water and air phases is modeled by deviation from Henry's law. The production of $gypsum$
    (weakened concrete) is modeled via a non-linear reaction rate $\eta$.
 Here we restrict our attention to a minimal set of chemical reactions
 mechanisms (as suggested in \cite{mbp}), namely

\begin{equation}\label{reeq}
 \left\{\begin{aligned}
10H^{+}+SO_4^{-2}+\mbox{org. matter} \;&\longrightarrow\;H_2S(aq)+4H_2O
+\mbox{oxid. matter}\\
H_2S(aq)+2O_2\;&{\longrightarrow}\;2H^++SO_4^{-2}\\
H_2S(aq)\;&\rightleftharpoons \;H_2S(g)\\
2H_2O+H^++SO_4^{-2}+CaCO_3\;&{\longrightarrow}\;CaSO_4\cdot 2H_2O+HCO_3^-
\end{aligned}
\right.
\end{equation}
We assume that reactions (\ref{reeq}) do not interfere with the mechanics of
the solid part of the pores. This is a rather strong assumption: We indicated already that (\ref{reeq}) can produce local
ruptures of the solid matrix \cite{taylor}, hence, generally we expect that the macroscopic
 mechanical properties of the piece of concrete will be affected. For more details on the
involved cement chemistry and practical aspects of acid corrosion, we
 refer  the reader to \cite{beddoe1} (for a nice enumeration of the involved physicochemical
 mechanisms),
\cite{taylor} (standard textbook on cement chemistry), as well as to
\cite{beddoe5,tixier} and references cited therein. For a
mathematical approach of a theme related to the conservation
and restoration of historical monuments [where the sulfatation reaction (\ref{reeq})
 plays an important role], we refer to the work by R.
Natalini and co-workers (cf. e.g. \cite{natalini}).
Based on the reactions mechanisms single-scale and two-scale (\ref{reeq}),
some models were derived in \cite{tasnim1,tasnim2,tasnim4,tasnim3}
mostly relying on periodic and locally-periodic homogenization. The model we
 investigate here has the same structure as one of these cases elucidated via
 homogenization.

\subsection{Geometry}
\begin{figure}[htbp]
\begin{center}
\includegraphics[width=.7\textwidth]{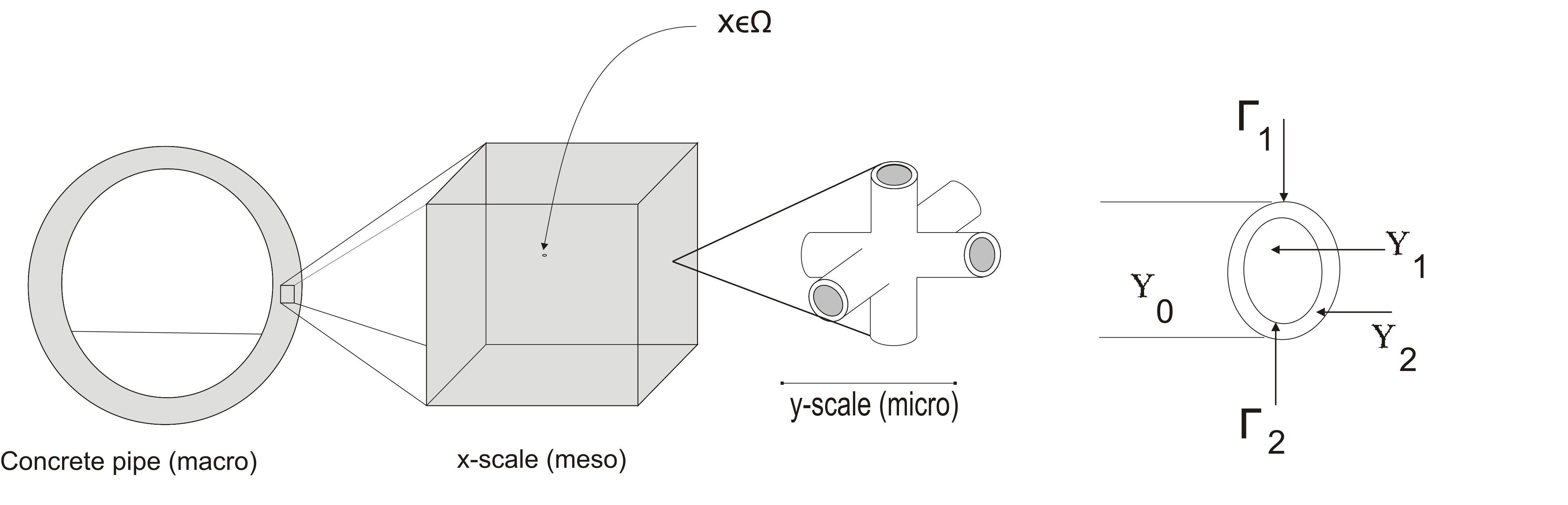}\\

\caption{\small\small{Left: Cross-section of a sewer pipe. Second from the left:
 A cubic piece from the concrete wall zoomed out. This is the scale we refer to as macroscopic.
 Second from the right: Reference pore configuration. Right: Zoomed one end of the cell.
}}\label{fig1}
\end{center}
\end{figure}

We consider $\Omega$ and $Y$ to be connected and bounded domains in $\mathds{R}^3$ defined on two different spatial scales.
$\Omega$ is a macroscopic domain with Lipschitz continuous boundary $\Gamma$ which is composed of two smooth disjoint parts: $\Gamma^N$ and $\Gamma^D$.
$Y$ is a standard unit cell associated with the microstructure within $\Omega$ with Lipschitz continuous boundary $\partial Y$ and has three disjoint components, i.e.,
 $Y:=\bar{Y}_0\cup \bar{Y}_1\cup \bar{Y}_2$, where $Y_0,{Y}_1,{Y}_2$ represent the solid matrix, the water layer which clings to the pipe wall
  and air-filled part surrounded by the water in the pipe, respectively.  All constituent parts of the pore connect neighboring
pores to one another.
 $\Gamma_1$, $\Gamma_2$ denote the inner boundaries within the cell $Y$, see Figure \ref{fig1},
that is, $\partial Y_1 = \Gamma_1 \cup \Gamma_2 \cup (\partial Y_1 \cap \partial Y)$.
$\Gamma_1$ represents the solid-water interface where the strong reaction takes place to destroy
the solid matrix and $\Gamma_2$ denotes the
water-air interface where the mass transfer occurs. $\Gamma_1, \Gamma_2$ are smooth enough
surfaces that do not touch each other.

\section{Two-scale model equations}\label{Model}

The two-scale reaction-diffusion system we have in mind consists of the following set of partial
differential equations coupled with one ordinary differential equation:
\begin{eqnarray}\label{main}
&& \partial_t w_1 - \nabla_y\cdot (d_1 \nabla_y w_1)= -f_1(w_1) +f_2(w_2) \quad
\text{ in }
(0,T)\times  \Omega \times Y_1, \label{pd1}\\
&&  \partial_t  w_2 - \nabla_y\cdot (d_2 \nabla_y w_2)= f_1(w_1) -f_2(w_2)
 \quad \text{ in } (0,T)\times\Omega \times Y_1, \label{pd2}\\
&& \partial_t  w_3 - \nabla\cdot (d_3 \nabla w_3)= -\alpha\int\limits_{\Gamma_2}\big(Hw_3-w_2\big)d\gamma_y
   \quad \text{ in } (0,T)\times\Omega, \label{pd3}\\
&&  \partial_t w_4 = \eta (w_1, w_4)  \quad
  \text{ on } (0,T)\times\Omega\times\Gamma_1.\label{od1}
\end{eqnarray}
 The system is equipped with the initial conditions
\begin{equation}\label{main_in}
\begin{cases}
 w_j(0, x,y) = w_j^0(x,y), \;\;j\in\{1,2\}  \hspace{1 cm} \text{ in } \Omega\times Y_1,  \\
w_3(0, x) = w_3^0(x)\quad  \text{ in } \Omega,  \qquad
 w_4(0, x,y) = w_4^0(x,y)  \hspace{0.4 cm} \text{ on } \Omega \times\Gamma_1,
\end{cases}
\end{equation}
while the boundary conditions are
\begin{equation}\label{main_bc}
\begin{cases}
d_1 \nabla_y w_1 \cdot \nu(y) =  -\eta (w_1, w_4) \quad  \text{ on } (0,T)\times\Omega \times\Gamma_1, \\
d_1 \nabla_y w_1 \cdot \nu(y) =  0  \quad \text{ on } (0,T)\times\Omega \times\Gamma_2 \mbox{ and }  (0,T)\times\Omega \times(\partial Y_1\cap\partial Y), \\
d_2 \nabla_y w_2 \cdot \nu(y) =  0 \quad \text{ on } (0,T)\times\Omega \times\Gamma_1 \mbox{ and }  (0,T)\times\Omega \times(\partial Y_1\cap\partial Y), \\
d_2 \nabla_y w_2 \cdot \nu(y) =\alpha(Hw_3-w_2\big)  \quad
 \text{ on }  (0,T)\times\Omega \times\Gamma_2,\\
d_3 \nabla w_3 \cdot \nu(x) =  0   \quad \text{ on } (0,T)\times\Gamma_N,\\
\;w_3 = w_3^D  \quad
\text{ on } (0,T)\times\Gamma_D,
\end{cases}
\end{equation}
where  $w_1$ denotes the concentration of $H_2SO_4$ in $(0,T)\times\Omega\times Y_1$, $w_2$ the concentration of $H_2S$ aqueous species in $(0,T)\times\Omega\times Y_1$,
$w_3$ the concentration of $H_2S$ gaseous species in $(0,T)\times\Omega$ and $w_4$ of $gypsum$ concentration on $(0,T)\times\Omega\times\Gamma_1$.
$\nabla$ without subscript denotes the differentiation w.r.t. macroscopic variable $x-$, while $\nabla_y, div_y$ are the respective differential operators w.r.t.  the micro-variable $y$. $\alpha$ denotes
the rate of the reaction taking place on the interface $\Gamma_2$
and $H$ is the Henry's constant. The microscale and macroscale are
 connected together via the right-hand side of $\eqref{main}_3$ and via the
micro-macro boundary condition $\eqref{main_bc}_4$. The information referring to the air phase $Y_2$
 is hidden in $w_3$.
The partial
differential equation for $w_3$, defined on macroscopic scale, is derived by averaging over $Y_2$; details are given in \cite{tasnim1}.

\section{Functional setting. Assumptions. Main results}\label{setting}

In this section, we enumerate the assumptions on the parameters and initial data needed to deal with the analysis of our problem. Furthermore, the definition
of the solution to the system \eqref{pd1}--\eqref{main_bc} is discussed and the main results of the paper are given at the end of this section.
To keep notation simple, we put $$X:=\{z\in H^1(\Omega)|z=0 \mbox{ on }\Gamma_D\}.$$
\begin{assumption}\label{assump}
\begin{enumerate}
\item[($A1$)] $d_i\in L^\infty(\Omega\times Y_1),i\in\{1,2\}$ and $d_3\in L^\infty(\Omega)$ such that \\
$(d_i(x,y)\xi,\xi)\geq d^0_{i}|\xi|^2$
  for $d^0_{i}>0$ for every $\xi\in\mathds{R}^3$, a.e. $(x,y)\in \Omega\times Y_1$ and $i\in\{1,2\}$, and
$(d_3(x)\xi,\xi)\geq d^0_{3}|\xi|^2$
  for $d^0_{3}>0$ for every $\xi\in\mathds{R}^3$ and a.e. $x \in \Omega$.

\item[($A2$)] $\eta(\alpha,\beta):= R(\alpha)Q(\beta)$, where $R$ and $Q$ are locally Lipschitz continuous functions such that
$R'\geq 0$ and $Q'\leq 0$ a.e. on $\mathds{R}$ and
   \begin{eqnarray}
R(\alpha):=
\left \{ \begin{array}{ccc}
\mbox{positive},\;\;\mbox{if}\;\;\alpha > 0,\\\nonumber
0,\;\;\;\mbox{otherwise, }\nonumber
 \end{array} \right.
\quad\quad
Q(\beta):=
\left \{ \begin{array}{ccc}
\mbox{positive},\;\;\mbox{if}\;\; \beta < \beta_{max},\\\nonumber
0,\;\;\;\mbox{otherwise, }\nonumber
 \end{array} \right.
\end{eqnarray}
where $\beta_{max}$ is a positive constant.
Also, we denote by $\hat{R}$ the primitive of $R$ with $\hat{R}(0) = 0$, that is,
$\hat{R} (r) = \int_0^r R(\xi) d\xi$ for $r \in \mathds{R}$.

\item[($A3$)]  The functions $f_i,i\in\{1,2\}$, are  increasing and locally Lipschitz continuous functions with $f_i(\alpha)= 0$ for $\alpha \leq 0$ and
$f_i(\alpha)> 0$ for $\alpha >0$, $i\in\{1,2\}$.
Furthermore, $\cR(f_1)=\cR(f_2)$, where $\cR(f)$ denotes the range of the function $f$.
Obviously, for  $M_1, M_2>0$ there exist positive constants $M'_1, M'_2 > 0$ such that
$$
f_1(M_1') = f_2(M'_2), M'_1 \geq M_1 \mbox{ and } M'_2 \geq M_2.
$$
\item[($A4$)] $w_{10}\in L^2(\Omega; H^1(Y_1))\cap L^\infty_+(\Omega\times Y_1),
w_{20} \in L^2(\Omega; H^1(Y_1))\cap L^\infty_+(\Omega\times Y_1)$, \\
$w_{30}\in H^1(\Omega)\cap L^\infty_+(\Omega)$, $w_{30} - w^D_3(0,\cdot) \in X$,
$w^D_{3}\in L^2(0,T;H^2(\Omega))\cap H^1(0, T;L^2(\Omega))\\\cap L_+^\infty((0,T)\times\Omega)$
with $\nabla w_3^D \cdot \nu = 0$ on $(0,T) \times \Gamma_N$,
 $w_{40}\in L^\infty_+(\Omega\times \Gamma_1)$.
\end{enumerate}
\end{assumption}

Note that in (A4) we define  $L_+^{\infty}(\Omega'): = L^{\infty}(\Omega') \cap \{u| u \geq 0 \mbox{ on } \Omega'\}$ for a domain $\Omega'$.
Next, we give the definition of a suitable concept of solution to our problem:

\begin{definition}\label{def}
We call the multiplet $(w_1, w_2, w_3, w_4)$ a solution to the problem \eqref{pd1}--\eqref{main_bc} if (S1) $\sim$ (S5) hold:
\begin{enumerate}
\item[(S1)]  $w_1, w_2 \in H^1(0,T; L^2(\Omega \times Y_1)) \cap L^{\infty}(0,T; L^2(\Omega; H^1(Y_1)))
\cap L^\infty((0,T)\times\Omega\times Y_1)$, \\
$w_3  \in H^1(0,T; L^2(\Omega)) \cap L^{\infty}((0,T) \times \Omega)$,
$w_3 - w_3^D \in L^\infty(0,T; X), \\
 w_4 \in H^1(0,T; L^2(\Omega\times\Gamma_1)) \cap L^\infty((0,T)\times\Omega\times \Gamma_1)$.

\item[(S2)] It holds that
\begin{eqnarray*}
&& \int_{\Omega \times Y_1} \partial_t w_1 ( w_1 - v_1) dxdy
+ \int_{\Omega \times Y_1} d_1 \nabla_y w_1 \cdot \nabla_y(w_1 - v_1) dxdy \\
& & + \int_{\Omega \times \Gamma_1} Q(w_4)(\hat{R}(w_1) - \hat{R}(v_1)) dxd\gamma_y \\
& \leq &
 \int_{\Omega \times Y_1} (- f_1(w_1) + f_2(w_2))  ( w_1 - v_1) dxdy \\
& & \qquad \qquad
\mbox{ for } v_1 \in L^2(\Omega; H^1(Y_1)) \mbox{ with } \hat{R}(v_1) \in L^1(\Omega \times \Gamma_1) \mbox{ a.e. on } [0,T].
\end{eqnarray*}

\item[(S3)] It holds that
\begin{eqnarray*}
&& \int_{\Omega \times Y_1} \partial_t w_2  v_2 dxdy
+ \int_{\Omega \times Y_1} d_2 \nabla_y w_2 \cdot \nabla_y v_2 dxdy
 - \alpha  \int_{\Omega \times \Gamma_2} (H w_3 - w_2) v_2 dxd\gamma_y \\
& = &
 \int_{\Omega \times Y_1} (f_1(w_1) - f_2(w_2))  v_2 dxdy
 \quad
\mbox{ for } v_2 \in L^2(\Omega; H^1(Y_1))  \mbox{ a.e. on } [0,T].
\end{eqnarray*}
\item[(S4)] It holds that
\begin{eqnarray*}
&& \int_{\Omega \times Y_1} \partial_t w_3  v_3 dx
+ \int_{\Omega } d_3 \nabla_y w_3 \cdot \nabla v_3 dxdy \\
& = & - \alpha  \int_{\Omega \times \Gamma_2} (H w_3 - w_2) v_3 dxd\gamma_y
 \quad
\mbox{ for } v_3 \in X \mbox{ a.e. on } [0,T].
\end{eqnarray*}

\item[(S5)] (\ref{od1}) holds a.e. on $(0,T) \times \Omega \times \Gamma_1$.
\end{enumerate}

\end{definition}
\begin{theorem}\label{Uniqueness}(Uniqueness)
Assume (A1)-(A4), then there exists at most one solution in
the sense of Definition \ref{def}.
\end{theorem}
{\it Proof.}
For the proof, see Section \ref{Proof}.
\hspace{.3 cm}
\begin{remark}
Having in view the proof of Theorem \ref{Uniqueness} and the working techniques in Theorem 3,
 pp. 520-521 in \cite{Dautray} as well as Theorem 4.1 in \cite{Muntean}, we expect that the solution
 in the sense of the Definition \ref{def} is
stable to the changes with respect to the initial data, boundary data, and model parameters.
\end{remark}


\begin{theorem}(Global existence of solutions to \eqref{pd1}--\eqref{main_bc})\label{existence}
Assume (A1)--(A4), then there exists a solution $(w_1,w_2, w_3, w_4)$ of the problem \eqref{pd1}--\eqref{main_bc}. Moreover, it holds that
\begin{itemize}
\item[(i)]$w_1(t),w_2(t)\geq0$ a.e. in ${\Omega\times Y_1}$, $w_3(t)\geq0$ a.e. in ${\Omega}$ and
 $w_4(t)\geq0$ a.e. on $\Omega\times\Gamma_1$ for a.e. $t \in [0,T]$.
\item[(ii)] $w_1(t)\leq M_1$, $w_2(t)\leq M_2$ a.e. in ${\Omega\times Y_1}$ ,
 $w_3(t)\leq M_3$ a.e. in ${\Omega}$ and
$w_4(t)\leq M_4$ a.e. on $\Omega\times\Gamma_1$ for a.e. $t \in [0,T]$,
 where $M_1$, $M_2$, $M_3$ and $M_4$ are positive constants satisfying
$M_1 \geq \|w_{10}\|_{L^{\infty}(\Omega \times Y_1)}$,
$M_2 \geq \|w_{20}\|_{L^{\infty}(\Omega \times Y_1)}$, \\
$M_3 \geq \max\{ \|w_{30}\|_{L^{\infty}(\Omega)}, \|w_{3}^D\|_{L^{\infty}(\Omega \times Y_1)}\}$,
$f_1(M_1) = f_2(M_2)$ and $M_2= HM_3$
and $M_4 = \max\{ \beta_{max}, \|w_{40}\|_{L^{\infty}(\Omega \times \Gamma_1)}\}$.
\end{itemize}
\end{theorem}

In order to prove the existence of a solution, we first solve the following problem P$_1(g, h)$ in Lemma
\ref{P1} for given functions $g$ and $h$:
\begin{eqnarray*}
&& \partial_t w_1 - \nabla_y\cdot (d_1 \nabla_y w_1) =  g \quad
\text{ in }
(0,T)\times  \Omega \times Y_1, \label{p1-1}\\
& & d_1 \nabla_y w_1 \cdot \nu(y) =  - hR(w_1) \quad  \text{ on } (0,T)\times\Omega \times\Gamma_1, \\
& & d_1 \nabla_y w_1 \cdot \nu(y) =  0  \quad \text{ on } (0,T)\times\Omega \times\Gamma_2 \mbox{ and }  (0,T)\times\Omega \times(\partial Y_1\cap\partial Y), \\
& & w_1(0) = w_{10} \quad \mbox{ on } \Omega \times Y_1.
\end{eqnarray*}

Next, for a given function $g$ on $(0,T) \times \Omega \times Y_1$,
we consider the following problem P$_2(g)$ (see Lemma \ref{P2}):
\begin{eqnarray*}
&& \partial_t w_1 - \nabla_y\cdot (d_1 \nabla_y w_1) =  g \quad
\text{ in }
(0,T)\times  \Omega \times Y_1, \label{p1-1}\\
& & d_1 \nabla_y w_1 \cdot \nu(y) =  - \eta(w_1, w_4) \quad  \text{ on } (0,T)\times\Omega \times\Gamma_1, \\
& & d_1 \nabla_y w_1 \cdot \nu(y) =  0  \quad \text{ on } (0,T)\times\Omega \times\Gamma_2 \mbox{ and }  (0,T)\times\Omega \times(\partial Y_1\cap\partial Y), \\
& & \partial_t w_4 = \eta(w_1, w_4) \quad \mbox{ a.e. on } (0,T) \times \Omega \times \Gamma_1, \\
& & w_1(0) = w_{10}  \mbox{ on } \Omega \times Y_1 \mbox{ and } w_{4}(0) = w_{40} \mbox{ on }
\Omega \times \Gamma_1.
\end{eqnarray*}

As a third step of the proof, we show the existence of a solution of the following problem P$_3(g)$ for a given function $g$ on $(0,T) \times \Omega \times Y_1$ (see Lemma \ref{P3}):
\begin{eqnarray*}
&&  \partial_t  w_2 - \nabla_y\cdot (d_2 \nabla_y w_2) = g
 \quad \text{ in } (0,T)\times\Omega \times Y_1, \label{p3-1}\\
&& \partial_t  w_3 - \nabla\cdot (d_3 \nabla w_3)=
 - \alpha \int\limits_{\Gamma_2} \big(Hw_3-w_2\big) d\gamma_y
   \quad \text{ in } (0,T)\times\Omega, \label{p3-2}\\
&& d_2 \nabla_y w_2 \cdot \nu =  0 \quad \text{ on } (0,T)\times\Omega \times\Gamma_1 \mbox{ and }  (0,T)\times\Omega \times(\partial Y_1\cap\partial Y), \\
&& d_2 \nabla_y w_2 \cdot \nu  = \alpha(Hw_3 - w_2)  \quad
 \text{ on }  (0,T)\times\Omega \times\Gamma_2,\\
&&d_3 \nabla w_3 \cdot \nu(x) =  0   \quad \text{ on } (0,T)\times\Gamma_N,\\
& & w_3 = w_3^D  \quad \text{ on } (0,T)\times\Gamma_D, \\
& & w_2(0) = w_{20}  \mbox{ on } \Omega \times Y_1 \mbox{ and } w_{3}(0) = w_{30} \mbox{ on }
   \Omega.
\end{eqnarray*}

\section{Auxiliary lemmas}

\begin{lemma} \label{P1}
Assume (A1), (A2), (A4),  $h \in H^1(0,T; L^2(\Omega \times \Gamma_1)) \cap L_+^{\infty}((0,T) \times \Omega \times \Gamma_1)$ and $g \in L^2((0,T) \times \Omega \times Y_1)$.
If $R$ is Lipschitz continuous and bounded on $\mathds{R}$,
 then there exists a solution $w_1$ of P$_1(g, h)$ in the following sense:
$w_1 \in H^1(0, T; L^2(\Omega \times Y_1)) \cap L^{\infty}(0, T; L^2(\Omega; H^1(Y_1)))$ satisfying \begin{eqnarray}
&& \int_{\Omega \times Y_1} \partial_t w_1 ( w_1 - v_1) dxdy
+ \int_{\Omega \times Y_1} d_1 \nabla_y w_1 \cdot \nabla_y(w_1 - v_1) dxdy \nonumber \\
& & + \int_{\Omega \times \Gamma_1} h (\hat{R}(w_1) - \hat{R}(v_1)) dxd\gamma_y \nonumber \\
& \leq &
 \int_{\Omega \times Y_1} g  ( w_1 - v_1) dxdy \quad
\mbox{ for } v_1 \in L^2(\Omega; H^1(Y_1))
\mbox{ a.e. on } [0,T],   \label{Lemma-1} \\
&& w_1(0) = w_{10} \quad \mbox{ on } \Omega \times Y_1.  \nonumber
\end{eqnarray}
\end{lemma}

{\it Proof. }
First, let  $\{\zeta_j\}$ be a Schauder basis of $L^2(\Omega; H^1(Y_1))$. More precisely,
$\{\zeta_j\}$ is an orthonormal system of a Hilbert space $L^2(\Omega \times Y_1)$ and is a fundamental of $L^2(\Omega; H^1(Y_1))$, that is, for any $z \in L^2(\Omega \times Y_1)$ we can take a sequence $\{z_k\}$ such that $z_k = \sum_{j=1}^{N_k} a_j^k \zeta_j$ and $z_k \to z$ in
$L^2(\Omega; H^1(Y_1))$ as $k \to \infty$, where $a_j^k \in \mathds{R}$.
Then there exists a sequence $\{w_{10}^n\}$ such that
$w_{10}^n := \sum_{j=1}^{N_n} \alpha_{j0}^n \zeta_j$ and $w_{10}^n \to w_{10}$ in $L^2(\Omega; H^1(Y_1))$ as $n \to \infty$.

Here, we are interested in the finite-dimensional approximations of the function $w_1$ that are of the form
\begin{equation}
w_1^n(t,x, y) := \sum_{j = 1}^{N_n} \alpha_j^n(t) \zeta(x,y)
\quad \mbox{ for } (t,x, y) \in (0,T) \times \Omega \times Y_1,
\label{fsol}
\end{equation}
where the coefficients $\alpha_j^n$, $j = 1, 2, \ldots, N_n$,  are determined by the following relations: For each $n$

\begin{eqnarray}
&& \int_{\Omega\times Y_1}
( \partial_t w^n_1(t) \phi_1  +
d_1 \nabla_y w^n_1(t) \nabla_y \phi_1 ) dxdy
+ \int_{\Omega \times \Gamma_1} h R(w_1^n(t)) \phi_1 dx d\gamma_y \nonumber\\
& = & \int_{\Omega \times Y_1} g(t) \phi_1  dxdy
\quad \mbox{ for } \phi_1 \in \mbox{ span}\{\zeta_{i}: i=1,2,..,N_n\} \mbox{ and }
 t \in (0,T], \label{ga1} \\
&& \alpha_j^n(0) = \alpha_{j0}^n \quad \mbox{ for } j = 1,2,...N_n. \label{ga2}
\end{eqnarray}
Consider $\phi_1 = \zeta_j$, $j = 1,2,...N_n$, as a test functions in (\ref{ga1}). This yields a system of ordinary differential equations
\begin{equation}
\partial_t \alpha^n_{j} (t) + \sum_{i=1}^{N_n} (A_{i})_{j} \alpha^n_{i} (t)
+  F_{j}^n (t, \alpha^n(t)) = J_j(t) \quad \mbox{ for } t \in (0,T] \mbox{ and } j = 1, 2, \ldots, N_n, \label{c1}
\end{equation}
where $\alpha^n(t) := (\alpha_1^n(t), \ldots, \alpha_{N_n}^n(t))$,
$(A_i)_j := \int_{\Omega \times Y_1} d_1 \nabla_y \zeta_i \cdot \nabla_y \zeta_j dxdy$, \\
$F_j^n(t, \alpha^n): = \int_{\Omega \times \Gamma_1} h(t) R(\sum_{i=1}^{N_n} \zeta_i) \zeta_j dx d\gamma_y$ and $J_j(t) = \int_{\Omega \times Y_1} g(t) \zeta_j dx dy$ for $t \in (0,T]$.
Note that $F_j^n$ is globally Lipschitz continuous due to the assumption of this lemma. According to the standard  existence theory for ordinary differential equations, there exists a unique  solution  $\alpha_j^n$,
$j = 1,2,..,N_n$,  satisfying \eqref{c1} for $0\leq t\leq T$ and \eqref{ga2}.
Thus the solution $w_1^n$ defined in \eqref{fsol} solves \eqref{ga1}--\eqref{ga2}.

Next, we show some uniform estimates for approximate solutions $w_1^n$ with respect to $n$.
We take $\phi_1 = w_1^n$  in \eqref{ga1} to obtain
\begin{eqnarray*}
 & &  \int_{\Omega\times Y_1}
 \partial_t w^n_1(t) w_1^n(t) dx dy  +  \int_{\Omega\times Y_1}
d_1 |\nabla_y w^n_1(t)|^2  dxdy
+ \int_{\Omega \times \Gamma_1} h(t) R(w_1^n(t)) w_1^n(t) dx d\gamma_y \nonumber\\
& = & \int_{\Omega \times Y_1} g(t) w_1^n(t)  dxdy
\quad \mbox{ for }  t \in (0,T]. \label{ga3}
\end{eqnarray*}
Since $R(r)r \geq 0$ for any $r \in \mathds{R}$, we see that
\begin{eqnarray*}
 & &  \frac{1}{2} \frac{d}{dt} \int_{\Omega\times Y_1}
  |w_1^n(t)|^2 dx dy  +  d_1^0\int_{\Omega\times Y_1} |\nabla_y w^n_1(t)|^2  dxdy
\nonumber\\
& \leq & \frac{1}{2} \int_{\Omega \times Y_1} |g(t)|^2 dx dy
+ \frac{1}{2} \int_{\Omega \times Y_1} |w_1^n(t)|^2  dxdy
\quad \mbox{ for }  t \in (0,T]. \label{ga3}
\end{eqnarray*}
Applying Gronwall's inequality, we have
$$
 \int_{\Omega\times Y_1}
  |w_1^n(t)|^2 dx dy  +  d_1^0 \int_0^t \int_{\Omega\times Y_1} |\nabla_y w^n_1(t)|^2  dxdy
 \leq C \quad \mbox{ for }  t \in (0,T],
$$
where $C$ is a positive constant independent of $n$.

To obtain  bounds on the time-derivative, we take $\phi_1 = \partial_t w_1^n$ as test function in \eqref{ga1}. It is easy to see that
\begin{eqnarray*}
 & &   \int_{\Omega\times Y_1} |\partial_t w_1^n(t)|^2 dx dy  +  \frac{1}{2} \frac{d}{dt}
 \int_{\Omega\times Y_1} d_1 |\nabla_y w^n_1(t)|^2  dxdy
+ \int_{\Omega \times \Gamma_1} h(t) \partial_t \hat{R}(w_1^n(t))  dx d\gamma_y
\nonumber\\
& = &  \int_{\Omega \times Y_1} g(t) \partial_t w_1^n(t) dx dy
  \\
& \leq  &  \frac{1}{2} \int_{\Omega \times Y_1} |g(t)|^2 dxdy
 + \frac{1}{2} \int_{\Omega \times Y_1} |\partial_t w_1^n(t)|^2 dx dy
\quad \mbox{ for a.e. }  t \in (0,T].
\end{eqnarray*}
Accordingly, we have
\begin{eqnarray*}
 & &   \frac{1}{2} \int_{\Omega\times Y_1} |\partial_t w_1^n|^2 dx dy
+  \frac{1}{2} \frac{d}{dt}
 \int_{\Omega\times Y_1} d_1 |\nabla_y w^n_1|^2  dxdy
+ \frac{d}{dt} \int_{\Omega \times \Gamma_1} h \hat{R}(w_1^n)  dx d\gamma_y
\nonumber\\
& \leq  &  \frac{1}{2} \int_{\Omega \times Y_1} |g|^2 dxdy
 +  \int_{\Omega \times \Gamma_1} \partial_t h \hat{R}(w_1^n) dx d\gamma_y
\quad \mbox{  a.e.  on }  \in (0,T].
\end{eqnarray*}
By integrating the latter equation, we have
\begin{eqnarray*}
 & &   \frac{1}{2} \int_0^{t_1} \int_{\Omega\times Y_1} \!\!|\partial_t w_1^n|^2 dx dy dt
+  \frac{d_1^0}{2} \int_{\Omega\times Y_1} \!\! |\nabla_y w^n_1(t_1)|^2  dxdy
+  \int_{\Omega \times \Gamma_1}\!\! h(t) \hat{R}(w_1^n(t_1))  dx d\gamma_y
\nonumber\\
& \leq  &    \frac{1}{2} \int_{\Omega\times Y_1}  d_1 |\nabla_y w^n_1(0)|^2  dxdy
+  \int_{\Omega \times \Gamma_1} h(0) \hat{R}(w_1^n(0))  dx d\gamma_y \\
&&+ \frac{1}{2} \int_0^{t_1} \int_{\Omega \times Y_1} |g|^2 dxdydt
 +  \int_0^{t_1} \int_{\Omega \times \Gamma_1} \partial_t h \hat{R}(w_1^n) dx d\gamma_ydt
\quad \mbox{ for }  t_1 \in (0,T].
\end{eqnarray*}
Note that
\begin{eqnarray*}
 \int_0^T \int_{\Omega \times \Gamma_1} |\hat{R}(w_1^n)|^2 dx d\gamma_y dt
& \leq &
  C \int_0^T \int_{\Omega \times \Gamma_1} |w_1^n|^2 dx d\gamma_y dt \\
& \leq &
  C \int_0^T \int_{\Omega \times Y_1} ( |\nabla_y w_1^n|^2 + |w_1^n|^2) dx dy dt.
\end{eqnarray*}
Here, we have used the trace inequality. Hence, we observe that
$\{w_1^n\}$ is bounded in $H^1(0,T; L^2(\Omega \times Y_1))$ and
$L^{\infty}(0,T; L^2(\Omega; H^1 ( Y_1))$. From these estimates we can choose a subsequence
$\{n_i\}$ of $\{n\}$ such that
$w_1^{n_i} \to w_1$ weakly in $H^1(0,T; L^2(\Omega \times Y_1))$, weakly* in
 $L^{\infty}(0,T; L^2(\Omega \times Y_1))$,  and  weakly* in $L^{\infty}(0, T; L^2(\Omega; H^1(Y_1)))$. Also, the above convergences implies that
$w_1^{n_i}(T) \to w_1(T)$ weakly in $L^2(\Omega \times Y_1)$.

Now, in order to show that \eqref{Lemma-1} holds let $v \in L^2(0,T; L^2(\Omega; H^1(Y_1)))$.
Obviously, we can take a sequence $\{v_k\}$ such that
$v_k(t) := \sum_{j=1}^{m_k} d_j^k(t) \zeta_j$ and $v_k \to v$ in
$L^2(0,T; L^2(\Omega; H^1(Y_1)))$ as $k \to \infty$, where $d_i^k \in C([0,T])$ for $i = 1, 2, \ldots, m_k$ and $k = 1, 2, \ldots$. For each $k$ and $i$ with $N_{n_i} \geq m_k$ from \eqref{ga1} it follows that
\begin{eqnarray*}
&& \int_0^T \int_{\Omega \times Y_1} \partial_t w_1^{n_i} ( w_1^{n_i} - v_1^k) dxdy dt
+ \int_{\Omega \times Y_1} d_1 \nabla_y w_1^{n_i} \cdot \nabla_y(w_1^{n_i} - v_1^k) dxdydt
       \nonumber \\
& & + \int_0^T
 \int_{\Omega \times \Gamma_1} h R(w_1^{n_i}) (w_1^{n_i} - v_1^k) dxd\gamma_ydt \nonumber \\
& \leq &
 \int_0^T \int_{\Omega \times Y_1} g  ( w_1^{n_i} - v_1^k) dxdydt.
\end{eqnarray*}
By the lower semi-continuity of the norm and the convex function $\hat{R}$ we have
\begin{eqnarray*}
&& \liminf_{i \to \infty} \left( \int_0^T \int_{\Omega \times Y_1} \{ \partial_t w_1^{n_i} ( w_1^{n_i} - v_1^k) +  d_1 \nabla_y w_1^{n_i} \cdot \nabla_y(w_1^{n_i} - v_1^k) \} dxdydt \right.
  \nonumber \\
& & \left.  \qquad \qquad +  \int_0^T
 \int_{\Omega \times \Gamma_1} h R(w_1^{n_i}) (w_1^{n_i} - v_1^k) dxd\gamma_ydt \right)
   \nonumber \\
& \geq &
  \int_0^T \int_{\Omega \times Y_1} \partial_t w_1 ( w_1 - v_1^k) dxdydt
+ \int_{\Omega \times Y_1} d_1 \nabla_y w_1 \cdot \nabla_y(w_1 - v_1^k) dxdydt \nonumber \\
& & + \int_0^T
 \int_{\Omega \times \Gamma_1} h (\hat{R}(w_1) - \hat{R}(v_1^k)) dxd\gamma_ydt  \nonumber \quad
 \mbox{ for each } k.
\end{eqnarray*}
Then we show that (\ref{Lemma-1}) holds for each $v_k$. Moreover, by letting $k \to \infty$
we obtain the conclusion of this lemma.
\hfill $\Box$

\vskip 12pt
Next, we solve the problem P$_2(g)$.
\begin{lemma} \label{P2}
Assume (A1), (A2), (A4),  and
$g \in L^2((0,T) \times \Omega \times Y_1)$.
If $R$ and $Q$ are Lipschitz continuous and bounded on $\mathds{R}$,
 then there exists a solution $(w_1, w_4)$ of P$_2(g)$ in the following sense:
$w_1 \in H^1(0, T; L^2(\Omega \times Y_1)) \cap L^{\infty}(0, T; L^2(\Omega; H^1(Y_1)))$
and $w_4 \in H^1(0,T; L^2(\Omega \times \Gamma_1))$
satisfying
\begin{eqnarray}
&& \int_{\Omega \times Y_1} \partial_t w_1 ( w_1 - v_1) dxdy
+ \int_{\Omega \times Y_1} d_1 \nabla_y w_1 \cdot \nabla_y(w_1 - v_1) dxdy \nonumber \\
& & + \int_{\Omega \times \Gamma_1} Q(w_4) (\hat{R}(w_1) - \hat{R}(v_1)) dxd\gamma_y \nonumber \\
& \leq &
 \int_{\Omega \times Y_1} g  ( w_1 - v_1) dxdy \quad
\mbox{ for } v_1 \in L^2(\Omega; H^1(Y_1))
\mbox{ a.e. on } [0,T],  \label{Lemma-2}  \\
&& \partial_t w_4 = \eta(w_1, w_4) \quad \mbox{ on }(0,T) \times \Omega \times \Gamma,
\label{lem2-s} \\
&& w_1(0) = w_{10} \quad \mbox{ on } \Omega \times Y_1  \mbox{ and }
w_4(0) = w_{40} \quad \mbox{ on } \Omega \times \Gamma_1. \nonumber
\end{eqnarray}
\end{lemma}
{\it Proof. }
Let $\bar{w}_4 \in  V := \{z \in  H^1(0,T; L^2(\Omega \times \Gamma_1)): z(0) = w_{40}\}$.
Then, since
$Q(\bar{w}_4) \in H^1(0,T; L^2(\Omega \times \Gamma_1)) \cap L_+^{\infty}((0,T) \times \Omega \times \Gamma_1)$, Lemma \ref{P1} implies that the problem P$_1(g, Q(\bar{w}_4))$ has a solution $w_1$ in the sense mentioned in Lemma \ref{P1}. Also, we put $w_4(t) := \int_0^t \eta(w_1(\tau), \bar{w}_4(\tau)) d\tau + w_{40}$
on $\Omega \times \Gamma_1$ for $t \in [0,T]$. Accordingly, we can define an operator
$\Lambda_T : V \to V$ by $\Lambda_T (\bar{w}_4) = w_4$.

Now, we  show that $\Lambda_{T}$ is a contraction mapping for sufficiently small $T > 0$.
Let $\bar{w}_4^i \in  H^1(0,T; L^2(\Omega \times \Gamma_1))$  and $w_1^i$ be a solution of
P$_1(g, Q(\bar{w}_4^i))$ and  $w_4^i = \Lambda_T(\bar{w}_4^i)$  for $i = 1, 2$, and
$w_1 = w_1^1 - w_1^2$, $w_4 = w_4^1 - w_4^2$ and  $\bar{w}_4 = \bar{w}_4^1 - \bar{w}_4^2$.

First, from \eqref{Lemma-1} with $v_1 = w_1^1$  we see that
\begin{eqnarray*}
&& \int_{\Omega \times Y_1} \partial_t w_1^1 ( w_1^1 - w_1^2) dxdy
+ \int_{\Omega \times Y_1} d_1 \nabla_y w_1^1 \cdot \nabla_y(w_1^1 - w_1^2) dxdy \nonumber \\
& & + \int_{\Omega \times \Gamma_1} Q(\bar{w}_4^1) (\hat{R}(w_1^1) - \hat{R}(w_1^2))
dxd\gamma_y \nonumber \\
& \leq &
 \int_{\Omega \times Y_1} g  ( w_1^1 - w_1^2) dxdy \quad
\mbox{ a.e. on } [0,T].    \label{l2-1}
\end{eqnarray*}
Similarly, we have
\begin{eqnarray*}
&& \int_{\Omega \times Y_1} \partial_t w_1^2 ( w_1^2 - w_1^1) dxdy
+ \int_{\Omega \times Y_1} d_1 \nabla_y w_1^2 \cdot \nabla_y(w_1^2 - w_1^1) dxdy \nonumber \\
& & + \int_{\Omega \times \Gamma_1} Q(\bar{w}_4^2) (\hat{R}(w_1^2) - \hat{R}(w_1^1))
dxd\gamma_y \nonumber \\
& \leq &
 \int_{\Omega \times Y_1} g  ( w_1^2 - w_1^2) dxdy \quad
\mbox{ a.e. on } [0,T].
\end{eqnarray*}
By adding these inequalities, for any $\varepsilon > 0$ we obtain
\begin{eqnarray}
&& \frac{1}{2} \frac{d}{dt} \int_{\Omega \times Y_1}  |w_1|^2  dxdy
+ \int_{\Omega \times Y_1} d_1 |\nabla_y w_1|^2  dxdy \nonumber \\
& \leq & -  \int_{\Omega \times \Gamma_1} (Q(\bar{w}_4^1) - Q(\bar{w}_4^2) )
      (\hat{R}(w_1^2) - \hat{R}(w_1^1)) dxd\gamma_y  \label{lem2-a} \\
& \leq & C_{\varepsilon}  \int_{\Omega \times \Gamma_1} |\bar{w}_4|^2  dxd\gamma_y
   + \varepsilon   \int_{\Omega \times \Gamma_1} |w_1|^2 dxd\gamma_y \nonumber  \\
& \leq & C_{\varepsilon}  \int_{\Omega \times \Gamma_1} |\bar{w}_4|^2  dxd\gamma_y
   + C_{Y_1} \varepsilon   \int_{\Omega \times Y_1} (|\nabla_y w_1|^2 + |w_1|^2) dxdy \quad
\mbox{ a.e. on } [0,T],   \nonumber
\end{eqnarray}
where $C_{Y_1}$ is a positive constant depending only on $Y_1$.
Here, by taking $\varepsilon > 0$ with $C_{Y_1} \varepsilon = \frac{1}{2}d_1^0$ and using Gronwall's inequality we see that
\begin{eqnarray}
& &  \frac{1}{2} \int_{\Omega \times Y_1}  |w_1(t)|^2 dxdy
 + \frac{d_1^0}{2}  \int_0^t \! \int_{\Omega \times Y_1} |\nabla_y w_1|^2 dxdy d\tau \nonumber \\
& \leq & e^{Ct}\! \int_0^t  \int_{\Omega \times \Gamma_1} \!|\bar{w}_4|^2 dx d\gamma_y d\tau \quad
\mbox{ for } t \in [0,T].  \label{lem2-3}
\end{eqnarray}

Next, on account of the definition of $w_4$ it is easy to see that
\begin{eqnarray*}
& & \frac{1}{2} \frac{d}{dt} \int_{\Omega \times \Gamma_1} |w_4(t)|^2 dx d\gamma_y \\
& \leq & \frac{1}{2}  \int_{\Omega \times \Gamma_1}
   (|\eta(w_1^1(t), \bar{w}_4^1(t)) - \eta(w_1^2(t), \bar{w}_4^2(t))|^2
+ |w_4(t)|^2)  dx d\gamma_y \\
& \leq & C  \int_{\Omega \times \Gamma_1}
   (|w_1(t)|^2 + |\bar{w}_4(t)|^2 + |w_4(t)|^2) dx d\gamma_y  \quad
  \mbox{ for a.e. } t \in [0, T].
\end{eqnarray*}
Gronwall's inequality, viewed in the context of \eqref{lem2-3},  implies that
\begin{eqnarray*}
& &   \int_{\Omega \times \Gamma_1} |w_4(t)|^2 dx d\gamma_y \\
& \leq &  C e^{Ct} (\int_0^t \int_{\Omega \times \Gamma_1}
   |\bar{w}_4|^2 dx d\gamma_y d\tau  +
\int_0^t \int_{\Omega \times Y_1}
   (|\nabla_y w_1|^2 + |w_1|^2) dx dy d\tau)
 \\
& \leq & C e^{Ct} \int_0^t \int_{\Omega \times \Gamma_1}
   |\bar{w}_4|^2 dx d\gamma_y d\tau  \quad
  \mbox{ for  } t \in [0, T].
\end{eqnarray*}
Hence, we obtain
\begin{eqnarray*}
  \|\partial_t {w}_4 \|_{L^2(0,T; L^2(\Omega \times \Gamma_1))}
& \leq &
\|\eta(w_1^1, \bar{w}_4^1) -  \eta(w_1^2, \bar{w}_4^2)\|_{L^2(0,T; L^2(\Omega \times \Gamma_1))}
\\
& \leq &  C (\|w_1\|_{L^2(0,T; L^2(\Omega; H^1(Y_1))} +
 \|\bar{w}_4\|_{L^2(0,T; L^2(\Omega  \times \Gamma_1))})
\\
& \leq &  C  \|\bar{w}_4\|_{L^2(0,T; L^2(\Omega  \times \Gamma_1))} \\
& \leq &  CT^{1/2}  \|\partial_t \bar{w}_4\|_{L^2(0,T; L^2(\Omega  \times \Gamma_1))},
\end{eqnarray*}
and
\begin{eqnarray*}
  \|\Lambda_T(\bar{w}_4^1) - \Lambda_T(\bar{w}_4^2)\|_{H^1(0,T; L^2(\Omega \times \Gamma_1))}
& \leq &     \|{w}_4 \|_{L^2(0,T; L^2(\Omega \times \Gamma_1))}
+ \|\partial_t {w}_4 \|_{L^2(0,T; L^2(\Omega \times \Gamma_1))}  \\
& \leq & C T^{1/2}  \|\bar{w}_4 \|_{H^1(0,T; L^2(\Omega \times \Gamma_1))}.
\end{eqnarray*}
This concludes that there exists $0 < T_0 \leq T$ such that $\Lambda_{T_0}$ is a contraction mapping. Here, we note that  the choice of $T_0$ is independent of initial values. Therefore, by applying Banach's fixed point theorem we have proved this lemma.
\hfill $\Box$

\vskip 12pt
As the third step of the proof of Theorem \ref{existence}, we solve P$_3(g)$.
\begin{lemma} \label{P3}
Assume (A1), (A2), (A4),  and
$g \in L^2((0,T) \times \Omega \times Y_1)$. Then there exists a pair $(w_2, w_3)$ such that
$w_2 \in H^1(0,T;  L^2(\Omega \times Y_1)) \cap L^{\infty}(0,T; L^2(\Omega; H^1(Y_1)))$,
$w_3 \in H^1(0,T;  L^2(\Omega)) \cap L^{\infty}(0,T; H^1(\Omega))$,
\begin{eqnarray}
&& \int_{\Omega \times Y_1} \partial_t w_2  v_2 dxdy
+ \int_{\Omega \times Y_1} d_2 \nabla_y w_2 \cdot \nabla_y v_2 dxdy
 - \alpha  \int_{\Omega \times \Gamma_2} (H w_3 - w_2) v_2 dxd\gamma_y \nonumber \\
& = &
 \int_{\Omega \times Y_1} g  v_2 dxdy
 \quad
\mbox{ for } v_2 \in L^2(\Omega; H^1(Y_1))  \mbox{ a.e. on } [0,T]. \label{lem3}
\end{eqnarray}
and  (S4).
\end{lemma}
{\it Proof. } Let $\{\zeta_j\}$ be the same set as in the proof of Lemma \ref{P1} and
$\{\mu_j\}$ be an orthonormal system of the Hilbert space $L^2(\Omega)$ and a fundamental of $X$.
Then we can take sequences $\{w_{20n}\}$ and $\{W_{30n}\}$ such that
$w_{20}^n := \sum_{j= 1}^{N_n} \beta_{j0}^n \zeta_j$,
$W_{30}^n := \sum_{j= 1}^{N_n} \beta_{j0}^n \mu_j$,
 $w_{20}^n \to w_{20}$ in $L^2(\Omega; H^1(Y_1))$ and
 $W_{30}^n \to w_{30} - w_3^D(0)$ in $X$ as $n \to \infty$.

We approximate $w_2$ and $W_3 := w_3 - w_3^D$ by functions $w_2^n$ and $W_3^n$ of the forms
\begin{equation}
w_2^n(t) = \sum_{j= 1}^{N_n} \beta_j^n(t) \zeta_j, \quad
 W_3^n(t) = \sum_{j= 1}^{N_n} \gamma_j^n(t) \mu_j \quad \mbox{ for } n, \label{lem3-6}
\end{equation}
where the coefficients $\beta_j^n$ and $\gamma_j^n$, $j = 1, 2, \ldots, N_n$ are determined the following relations: For each $n$ , we have:
\begin{eqnarray}
&& \int_{\Omega \times Y_1} \partial_t w_2^n  \phi_2 dxdy
+ \int_{\Omega \times Y_1} d_2 \nabla_y w_2^n \cdot \nabla_y \phi_2 dxdy
 - \alpha  \int_{\Omega \times \Gamma_2} (H W_3^n - w_2^n) \phi_2 dxd\gamma_y \nonumber \\
& = &
 \int_{\Omega \times Y_1} g  \phi_2 dxdy
+ \alpha \int_{\Omega \times \Gamma_2}\!\! H w_3^D  \phi_2 dxd\gamma_y  \label{lem3-1} \\
& & \quad
\qquad \quad \mbox{ for } \phi_2 \in \mbox{ span}\{\zeta_{i}: i=1,.., N_n\},  t \in (0,T],
     \nonumber \\
&& \beta_j^n(0) = \beta_{j0}^n,  \nonumber
\end{eqnarray}
\begin{eqnarray}
&& \int_{\Omega \times Y_1} \partial_t W_3^n  \phi_3 dx
+ \int_{\Omega } d_3 \nabla W_3^n \cdot \nabla \phi_3 dx
+ \alpha  \int_{\Omega \times \Gamma_2} (H W_3^n - w_2^n) \phi_3 dxd\gamma_y \nonumber \\
& = &
  -\int_{\Omega} (\partial_t w_3^D - \nabla d_3 (\nabla w_3^D)) \phi_3 dx dy
- \alpha  \int_{\Omega \times \Gamma_2} H w_3^D  \phi_3 dxd\gamma_y \label{lem3-2} \\
& &   \quad
\qquad \quad \mbox{ for } \phi_3 \in \mbox{ span}\{\mu_{i}: i=1,2,..., N_n\},
 t \in (0,T], \nonumber  \\
&& \gamma_j^n(0) = \gamma_{j0}^n \quad \mbox{ for } j = 1,2,..,N_n. \nonumber
\end{eqnarray}
Consider $\phi_2 = \zeta_j$ and $\phi_3 = \mu_j$, $j = 1,2,..,N_n$, as a test functions in
(\ref{lem3-1}) and \eqref{lem3-2}, respectively, these yield a system of ordinary differential equations
\begin{eqnarray*}
& & \partial_t \beta^n_{j} (t) + \sum_{i=1}^{N_n} (B_{i})_{j} \beta^n_{i} (t)
+  \sum_{i = 1}^{N_n} (\tilde{B}_i)_j  \gamma_i^n(t) = J_{j2}(t)  \mbox{ for } t \in (0,T]
  \mbox{ and } j = 1, 2, \ldots, N_n, \label{lem3-3} \\
& &  \partial_t \gamma^n_{j} (t) + \sum_{i=1}^{N_n} (C_{i})_{j} \gamma^n_{i} (t)
+  \sum_{i = 1}^{N_n} (\tilde{C}_i)_j  \beta_i^n(t) = J_{j3}(t)  \mbox{ for } t \in (0,T]
  \mbox{ and } j = 1, 2, \ldots, N_n, \label{lem3-3}
\end{eqnarray*}
where
$(B_i)_j := \int_{\Omega \times Y_1} d_2 \nabla_y \zeta_i \cdot \nabla_y \zeta_j dxdy$,
$(\tilde{B}_i)_j := \int_{\Omega \times \Gamma_2} \mu_i \zeta_j dxd\gamma_y$, \\
$(C_i)_j := \int_{\Omega} d_3 \nabla \mu_i \cdot \nabla \mu_j dx +
 \alpha H \int_{\Omega \times \Gamma_2} \mu_i \mu_j dx d\gamma_y$,
$(\tilde{C}_i)_j := - \alpha \int_{\Omega \times \Gamma_2} \mu_j \zeta_i dxd\gamma_y$, \\
$J_{2j}(t): =  \int_{\Omega \times Y_1} g(t) \zeta_j dx dy$,
$J_{3j}(t): =  \int_{\Omega \times Y_1} (\partial_t w_3^D - \nabla (\nabla w_3^D)) \zeta_j dx dy
+ \alpha H \int_{\Omega \times \Gamma_2} w_3^D \zeta_j dx d\gamma_y$ for $t \in (0,T]$.

Clearly, this linear system of ordinary differential equations has a solution $\beta_j^n$ and
$\gamma_j^n$. Thus the solutions $w_2^n$ and $W_3^n$ defined in \eqref{lem3-6} solve
(\ref{lem3-1}) and \eqref{lem3-2},  respectively.

Next, we shall obtain some uniform estimates for $w_2^n$ and $W_3^N$. We take $\phi_2 = w_2^n$ and
$\phi_3 = W_3^n$ in \eqref{lem3-1} and \eqref{lem3-2}, respectively, to have
\begin{eqnarray*}
&& \frac{1}{2} \frac{d}{dt} \int_{\Omega \times Y_1} |w_2^n(t)|^2 dxdy
 +  \int_{\Omega \times Y_1}  d_2 |\nabla_y w_2^n(t)|^2 dxdy
 + \frac{\alpha}{4}  \int_{\Omega \times \Gamma_2}  |w_2^n(t)|^2 dxd\gamma_y\\
& \leq &
  \frac{\alpha}{2}H |\Gamma_2| \int_{\Omega} |W_3^n(t)|^2 dx
 +  \frac{1}{2} \int_{\Omega \times Y_1} |g(t)|^2 dxdy   \\
& &  +  \frac{1}{2}  \int_{\Omega \times Y_1} |w_2^n(t)|^2 dxdy
 +  \alpha H^2  \int_{\Omega \times \Gamma_2} |w_3^D(t)|^2 dxdy \quad \mbox{ for a.e. } t \in [0,T],\end{eqnarray*}
\begin{eqnarray*}
&& \frac{1}{2} \frac{d}{dt} \int_{\Omega} |W_3^n(t)|^2 dx
 +  \int_{\Omega}  d_3 |\nabla W_3^n(t)|^2 dx
 + \alpha H  |\Gamma_2| \int_{\Omega}  |W_3^n(t)|^2 dxd\gamma_y\\
& \leq &
  \frac{1}{4} \int_{\Omega \times \Gamma_2} |w_2^n(t)|^2 dxd\gamma_y
 +  (\alpha^2 |\Gamma_2| + \frac{1}{2})  \int_{\Omega} |W_3^n(t)|^2 dx \\
&& +  \frac{1}{2}  \int_{\Omega} |g^D(t)|^2 dx \quad  \mbox{ for a.e. } t \in [0,T],
\end{eqnarray*}
where $|\Gamma_2| := \int_{\Gamma_2} d\gamma_y$, and
$g^D := \partial_t w_3^D - \nabla (\nabla w_3^D) + \alpha H |\Gamma_2| w_3^D$.
By adding these inequalities, we get
\begin{eqnarray*}
&& \frac{1}{2} \frac{d}{dt} \int_{\Omega \times Y_1} |w_2^n(t)|^2 dxdy
 + d_2^0 \int_{\Omega \times Y_1}   |\nabla_y w_2^n(t)|^2 dxdy \\
& & +  \frac{1}{2} \frac{d}{dt} \int_{\Omega} |W_3^n(t)|^2 dx
 + d_3^0 \int_{\Omega}  |\nabla W_3^n(t)|^2 dx  \\
& \leq &
(\alpha^2 |\Gamma_2| + \frac{1}{2} + \frac{H}{2} |\Gamma_2|)  \int_{\Omega} |W_3^n(t)|^2 dx
 +  \frac{1}{2} \int_{\Omega \times Y_1} |g(t)|^2 dxdy   \\
&&  +  \frac{3}{4}  \int_{\Omega \times Y_1} |w_2^n(t)|^2 dxdy
 +  \frac{1}{2}  \int_{\Omega} |g^D(t)|^2 dx \quad  \mbox{ for a.e. } t \in [0,T].
\end{eqnarray*}
Consequently, Gronwall's inequality implies that for some positive constant $C$
\begin{eqnarray}
&&  \int_{\Omega \times Y_1} |w_2^n(t)|^2 dxdy
 +  \int_{\Omega} |W_3^n(t)|^2 dx \leq C  \quad  \mbox{ for } t \in [0,T] \mbox{ and } n,
  \label{lem3-9} \\
&&
  \int_0^T \int_{\Omega \times Y_1}   |\nabla_y w_2^n|^2 dxdydt
 + \int_0^T \int_{\Omega}  |\nabla W_3^n(t)|^2 dx dt \leq C \mbox{ for } n. \label{lem3-10}
\end{eqnarray}

To obtain a uniform estimate for the time derivative by taking $\phi_2 = \partial_t  w_2^n$
in \eqref{lem3-1}  we observe that
 \begin{eqnarray}
&&  \int_{\Omega \times Y_1} |\partial_t w_2^n(t)|^2 dxdy
 + \frac{1}{2} \frac{d}{dt} \int_{\Omega \times Y_1}  d_2 |\nabla_y w_2^n(t)|^2 dxdy \nonumber \\
& & - \int_{\Omega \times \Gamma_2}  (H(W_3^n(t) + w_3^D(t))- w_2^n(t)) \partial_t w_2^n(t) dxd\gamma_y \nonumber
\\
& =  &
  \int_{\Omega \times Y_1}  g(t) \partial_t w_2^n(t) dxdy
 \quad  \mbox{ for a.e. } t \in [0,T]. \label{xx}
\end{eqnarray}
Here, we denote the third term in the left hand side of \eqref{xx}  by $J(t)$ and see that \begin{eqnarray*}
J(t) & := &
- \frac{d}{dt} \int_{\Omega \times \Gamma_2} H(W_3^n(t) + w_3^D(t)) w_2^n(t) dxd\gamma_y +
\frac{1}{2} \frac{d}{dt} \int_{\Omega \times \Gamma_2} |w_2^n(t)|^2 dxd\gamma_y \\
& & + H \int_{\Omega \times \Gamma_2}
  \partial_t (W_3^n(t) + w_3^D(t)) w_2^n(t) dxd\gamma_y  \quad  \mbox{ for a.e. } t \in [0,T].
\end{eqnarray*}
Then we have
 \begin{eqnarray*}
&&  \frac{1}{2} \int_{\Omega \times Y_1} |\partial_t w_2^n(t)|^2 dxdy
 + \frac{1}{2} \frac{d}{dt} \int_{\Omega \times Y_1}  d_2 |\nabla_y w_2^n(t)|^2 dxdy +
\frac{1}{2} \frac{d}{dt} \int_{\Omega \times \Gamma_2} |w_2^n(t)|^2 dxd\gamma_y \\
& \leq  & \frac{1}{2} \int_{\Omega \times Y_1}  |g(t)|^2 dxdy
+  \frac{d}{dt} \int_{\Omega \times \Gamma_2} H(W_3^n(t) + w_3^D(t)) w_2^n(t) dxd\gamma_y \\
& &  + H |\Gamma_2|^{1/2} \int_{\Omega}
  |\partial_t (W_3^n(t) + w_3^D(t))|  \|w_2^n(t)\|_{L^2(\Gamma_2)} dx
 \quad  \mbox{ for a.e. } t \in [0,T].
\end{eqnarray*}
Similarly, by taking $\phi_3 = \partial_t  W_3^n$
in \eqref{lem3-2}  we have
 \begin{eqnarray*}
&&  \int_{\Omega} |\partial_t W_3^n(t)|^2 dx
 + \frac{1}{2} \frac{d}{dt} \int_{\Omega}  d_3 |\nabla W_3^n(t)|^2 dx \\
&=  & -  \alpha \int_{\Omega \times \Gamma_2}  (H(W_3^n(t) + w_3^D(t)) - w_2^n(t)) \partial_t W_3^n(t) dxd\gamma_y
-   \int_{\Omega}  g^D(t) \partial_t W_3^n(t) dx   \\
& \leq & -  \alpha H |\Gamma_2| \frac{d}{dt} \int_{\Omega \times \Gamma_2} |W_3^n(t)|^2 dx
+ 2\alpha^2 H^2 |\Gamma_2|^2 \int_{\Omega} |w_3^D(t)|^2 dx \\
& & + 2 \alpha^2 |\Gamma_2|  \int_{\Omega \times \Gamma_2} |w_2^n(t)|^2 dx d\gamma_y
+ \frac{1}{2} \int_{\Omega} |\partial_t W_3^n(t)|^2 dx
+  \int_{\Omega}  |g^D(t)|^2 dx
  \mbox{ for a.e. } t \in [0,T].
\end{eqnarray*}
From these inequalities,  it follows that
 \begin{eqnarray*}
&&  \frac{1}{2} \int_{\Omega \times Y_1} |\partial_t w_2^n(t)|^2 dxdy
 + \frac{1}{2} \frac{d}{dt} \int_{\Omega \times Y_1}  d_2 |\nabla_y w_2^n(t)|^2 dxdy
+ \frac{1}{2} \frac{d}{dt} \int_{\Omega \times \Gamma_2} |w_2^n(t)|^2 dxd\gamma_y \\
&&  +  \frac{1}{8} \int_{\Omega} |\partial_t W_3^n(t)|^2 dx
 + \frac{1}{2} \frac{d}{dt} \int_{\Omega}  d_3 |\nabla W_3^n(t)|^2 dx
+ \alpha H |\Gamma_2| \frac{d}{dt} \int_{\Omega \times \Gamma_2} |W_3^n(t)|^2 dx \\
& \leq &
\frac{1}{2} \int_{\Omega \times Y_1}  |g(t)|^2 dxdy
+  \frac{d}{dt} \int_{\Omega \times \Gamma_2} H(W_3^n(t) + w_3^D(t)) w_2^n(t) dxd\gamma_y \\
& & + (2 |\Gamma_2| +  2 \alpha^2 |\Gamma_2|)
   \int_{\Omega \times \Gamma_2} |w_2^n(t)|^2 dx d\gamma_y
+  \int_{\Omega}  |g^D(t)|^2 dx   \quad
  \mbox{ for a.e. } t \in [0,T].
\end{eqnarray*}
Here, we use Gronwall's inequality, again, and have
 \begin{eqnarray*}
&&  \frac{1}{2} \int_0^{t_1} \int_{\Omega \times Y_1} |\partial_t w_2^n|^2 dxdydt
 + \frac{1}{2}  \int_{\Omega \times Y_1}  d_2 |\nabla_y w_2^n(t_1)|^2 dxdy
+ \frac{1}{2}  \int_{\Omega \times \Gamma_2} |w_2^n(t_1)|^2 dxd\gamma_y \\
&&  +  \frac{1}{8} \int_0^{t_1} \int_{\Omega} |\partial_t W_3^n|^2 dxdt
 + \frac{1}{2}  \int_{\Omega}  d_3 |\nabla W_3^n(t_1)|^2 dx
+ \alpha H |\Gamma_2|  \int_{\Omega \times \Gamma_2} |W_3^n(t_1)|^2 dx \\
& \leq & e^{Ct_1} \int_0^{t_1}\left( \int_{\Omega \times Y_1}  |g|^2 dxdy
  + \int_{\Omega}  |g^D|^2 dx\right) dt  \\
& & + e^{Ct_1}
\int_0^{t_1} e^{- Ct}( \frac{d}{dt} \int_{\Omega \times \Gamma_2} H(W_3^n + w_3^D) w_2^n
  dxd\gamma_y ) dt \\
& \leq & e^{Ct_1} \int_0^{t_1}(\int_{\Omega \times Y_1}  |g|^2 dxdy
  + \int_{\Omega}  |g^D|^2 dx) dt
  +  \int_{\Omega \times \Gamma_2} H(W_3^n(t_1) + w_3^D(t_1)) w_2^n(t_1) dxd\gamma_y \\
& & + e^{Ct_1} |\int_{\Omega \times \Gamma_2} H(W_3^n(0) + w_3^D(0)) w_2^n(0) dxd\gamma_y| \\
& & + e^{Ct_1} \int_0^{t_1} \int_{\Omega \times \Gamma_2} H(W_3^n + w_3^D) w_2^n dxd\gamma_y dt \quad \mbox{ for } t_1 \in [0, T].
 \end{eqnarray*}
This inequality together with \eqref{lem3-9} and \eqref{lem3-10} leads to
\begin{eqnarray}
&&     \int_{\Omega \times Y_1}   |\nabla_y w_2^n(t) |^2 dxdy +
\int_{\Omega}  |\nabla W_3^n(t)|^2 dx
 \leq C  \quad  \mbox{ for } t \in [0,T] \mbox{ and } n,
  \label{lem3-11} \\
&&
  \int_0^T \int_{\Omega \times Y_1}   |\partial_t w_2^n|^2 dxdydt
 + \int_0^T \int_{\Omega}  |\partial_t W_3^n|^2 dx dt \leq C \mbox{ for } n. \label{lem3-12}
\end{eqnarray}

By \eqref{lem3-9} $\sim$ \eqref{lem3-12} there exists a subsequence $\{n_i\}$  such that
$w_2^{n_i} \to w_2$ weakly in $H^1(0,T; L^2(\Omega \times Y_1))$, weakly* in
$L^{\infty}(0,T; L^2(\Omega; H^1(Y_1))$ and
$W_3^{n_i} \to W_3$ weakly in \\$H^1(0,T; L^2(\Omega))$, weakly* in $L^{\infty}(0,T; H^1(\Omega))$
as $i \to  \infty$. Clearly,
$w_2^{n_i} \to w_2$ weakly in $L^2((0,T) \times \Omega \times \Gamma_2)$ as $i \to \infty$.
Here, we put $w_3 = W_3 + w_3^D$.

Since the problem P$_3(g)$ is linear, similarly to the last part of the proof of Lemma \ref{P1},
we can show \eqref{lem3} and (S4). \hfill $\Box$

\section{Proof of our main results} \label{Proof}
First, we consider our problem \eqref{pd1}--\eqref{main_bc} in the case when $f_1$, $f_2$, $R$ and
$Q$ are Lipschitz continuous and bounded on  $\mathds{R}$.

\begin{proposition} \label{Lipschitz}
If (A1)-(A4) hold and $f_1$, $f_2$, $R$ and $Q$ are Lipschitz continuous and bounded on
$\mathds{R}$,
then there exists one and only one multiplet $(w_1, w_2, w_3, w_4)$ satisfying
$$ \left\{ \begin{array}{l}
w_1, w_2 \in H^1(0,T; L^2(\Omega \times Y_1)) \cap L^{\infty}(0,T; L^2(\Omega; H^1(Y_1))), \\
w_3  \in H^1(0,T; L^2(\Omega)), w_3 - w_3^D \in L^\infty((0,T; X),
w_4 \in H^1(0,T; L^2(\Omega\times\Gamma_1)),  \\
\mbox{ (S2) holds for any } v_1 \in L^2(\Omega; H^1(Y_1)),  \mbox{ and (S3), (S4) and (S5) hold.}
\end{array}
\right.
\leqno{\mbox{(S')}}
$$
 \end{proposition}
{\it Proof. } Let $(\bar{w}_1, \bar{w}_2) \in L^2((0,T) \times \Omega \times Y_1)^2$.
Then,  by Lemmas \ref{P2} and \ref{P3},  there exist solutions $(w_1, w_4)$ of P$_2(-f_1(\bar{w}_1) + f_2(\bar{w}_2))$ and $(w_2, w_3)$ of P$_3(f_1(\bar{w}_1) - f_2(\bar{w}_2))$, respectively.
Accordingly, we can define an operator $\bar{\Lambda}_{T}$ from
$L^2((0,T) \times \Omega \times Y_1)^2$ into itself. From now on, we show that $\bar{\Lambda}_T$ is contraction for small $T$. To do so, let $(\bar{w}_1^i, \bar{w}_2^i) \in L^2((0,T) \times \Omega \times Y_1)^2$,   $(w_1^i, w_4^i)$ and $(w_2^i, w_3^i)$ be solutions
 of P$_2(-f_1(\bar{w}_1^i) + f_2(\bar{w}_2^i))$ and P$_3(f_1(\bar{w}_1^i) - f_2(\bar{w}_2^i))$, respectively, for $i = 1, 2$, and put
$\bar{w}_1 = \bar{w}_1^1 - \bar{w}_1^2$,
$\bar{w}_2 = \bar{w}_2^1 - \bar{w}_2^2$,
$w_j = w_j^1 - w_j^2$, $j = 1, 2, 3, 4$.

Similarly to \eqref{lem2-a}, we see that
\begin{eqnarray}
&& \frac{1}{2} \frac{d}{dt} \int_{\Omega \times Y_1}  |w_1|^2  dxdy
+ \int_{\Omega \times Y_1} d_1 |\nabla_y w_1|^2  dxdy \nonumber \\
& \leq & -  \int_{\Omega \times \Gamma_1} (Q(\bar{w}_4^1) - Q(\bar{w}_4^2) )
      (\hat{R}(w_1^2) - \hat{R}(w_1^1)) dxd\gamma_y  \nonumber\\
&  & -  \int_{\Omega \times Y_1} (f_1(\bar{w}_1^1) - f_1(\bar{w}_1^2)) w_1 dxdy
     + \int_{\Omega \times Y_1} (f_2(\bar{w}_1^1) - f_2(\bar{w}_1^2)) w_1 dxdy \nonumber \\
& \leq &
\frac{d_1^0}{2} \int_{\Omega \times Y_1} |\nabla_y w_1|^2  dxdy
+ C \int_{\Omega \times Y_1} |w_1|^2  dxdy
+ C \int_{\Omega \times \Gamma_1} |w_4|^2  dxd\gamma_y  \nonumber \\
& & + C \int_{\Omega \times Y_1} (|\bar{w}_1|^2 + |\bar{w}_2|^2)  dxdy
\quad
\mbox{ a.e. on } [0,T].   \nonumber
\end{eqnarray}

Next, we test \eqref{lem3} by $w_2$. Consequently, by elementary calculations, we obtain
\begin{eqnarray*}
& & \frac{1}{2} \frac{d}{dt} \int_{\Omega \times Y_1} |w_2|^2 dxdy
+ d_2^0 \int_{\Omega \times Y_1} |\nabla_y w_2|^2 dxdy
 + \alpha \int_{\Omega \times \Gamma_2} |w_2|^2 dxd\gamma_y \\
& \leq &
\int_{\Omega \times Y_1} (f_1(\bar{w}_1^1) - f_1(\bar{w}_1^2)) w_2 dxdy
- \int_{\Omega \times Y_1} (f_1(\bar{w}_2^1) - f_1(\bar{w}_2^2)) w_2 dxdy \\
& & + \alpha \int_{\Omega \times \Gamma_2} Hw_3 w_2 dx d\gamma_y \\
& \leq &
C \int_{\Omega \times Y_1} (|\bar{w}_1| + |\bar{w}_2|) |w_2| dxdy
 + \frac{\alpha}{2} \int_{\Omega \times \Gamma_2} |w_2|^2 dx d\gamma_y
 + \frac{\alpha}{2}H^2 |\Gamma_2| \int_{\Omega} |w_3|^2 dx
\end{eqnarray*}
and
\begin{eqnarray*}
& & \frac{1}{2} \frac{d}{dt} \int_{\Omega \times Y_1} |w_2|^2 dxdy
+ d_2^0 \int_{\Omega \times Y_1} |\nabla_y w_2|^2 dxdy
 + \frac{\alpha}{2} \int_{\Omega \times \Gamma_2} |w_2|^2 dxd\gamma_y \\
& \leq &
C \int_{\Omega \times Y_1} (|\bar{w}_1|^2 + |\bar{w}_2|^2 + |w_2|^2) dxdy
 + C \int_{\Omega} |w_3|^2 dx  \quad \mbox{ a.e. on } [0,T].
\end{eqnarray*}

It follows form (S4) that
\begin{eqnarray*}
& & \frac{1}{2} \frac{d}{dt} \int_{\Omega} |w_3|^2 dx
+ d_3^0 \int_{\Omega} |\nabla w_3|^2 dx
+ \alpha H \int_{\Omega \times \Gamma_2} |w_3|^2 dxd\gamma_y
\\
& \leq &
\frac{\alpha}{4}  \int_{\Omega \times \Gamma_2} |w_2|^2 dxd\gamma_y
 + \alpha |\Gamma_2| \int_{\Omega} |w_3|^2 dx  \quad \mbox{ a.e. on } [0,T].
\end{eqnarray*}

Moreover, by using the trace inequality and \eqref{lem2-s},
we see that for $\varepsilon > 0$   we can write
\begin{eqnarray*}
& & \frac{1}{2} \frac{d}{dt} \int_{\Omega \times \Gamma_1} |w_4|^2 dx \\
& \leq &
  \int_{\Omega \times \Gamma_1} |\eta(w_1^1, w_4^1) - \eta(w_1^2, w_4^2)| |w_4|  dxd\gamma_y \\
& \leq &
 C  \int_{\Omega \times \Gamma_1} (|w_1||w_4| + |w_4|^2)  dxd\gamma_y \\
& \leq &
C_{Y_1} \varepsilon \int_{\Omega \times Y_1}(|\nabla_y w_1|^2 + |w_1|^2) dxdy
+  C  \int_{\Omega \times \Gamma_1}  |w_4|^2  dxd\gamma_y
\quad \mbox{ a.e. on } [0,T].
\end{eqnarray*}

Here, we take $\varepsilon$ with $C_{Y_1} \varepsilon = \frac{d_1^0}{4}$ and add the above inequalities. Then it holds that
\begin{eqnarray*}
&& \frac{1}{2} \frac{d}{dt} \int_{\Omega \times Y_1}  |w_1|^2  dxdy
+ \frac{d_1^0}{4} \int_{\Omega \times Y_1}  |\nabla_y w_1|^2  dxdy \nonumber \\
& & + \frac{1}{2} \frac{d}{dt} \int_{\Omega \times Y_1} |w_2|^2 dxdy
+ d_2^0 \int_{\Omega \times Y_1} |\nabla_y w_2|^2 dxdy
 + \frac{\alpha}{4} \int_{\Omega \times \Gamma_2} |w_2|^2 dxd\gamma_y \\
& & + \frac{1}{2} \frac{d}{dt} \int_{\Omega} |w_3|^2 dx
+ d_3^0 \int_{\Omega} |\nabla w_3|^2 dx
+ \alpha H \int_{\Omega \times \Gamma_2} |w_3|^2 dxd\gamma_y
 + \frac{1}{2} \frac{d}{dt} \int_{\Omega \times \Gamma_1} |w_4|^2 dx \\
& \leq &
C \int_{\Omega \times Y_1} (|\bar{w}_1|^2 + |\bar{w}_2|^2) dx dy
 + C \int_{\Omega \times Y_1} (|w_1|^2 + |w_2|^2) dxdy \\
& &  + C \int_{\Omega} |w_3|^2 dx
 + C \int_{\Omega \times \Gamma_1} |w_4|^2 dxdy
 \quad \mbox{ a.e. on } [0,T].
\end{eqnarray*}
Hence, Gronwall's inequality implies that
\begin{eqnarray*}
& &  \int_{\Omega \times Y_1}  (|w_1(t)|^2 + |w_2(t)|^2) dxdy
+ \int_{\Omega} |w_3(t)|^2 dx
+ \int_{\Omega \times \Gamma_1} |w_4(t)|^2 dxd\gamma_y \\
& \leq &
e^{Ct} \int_0^t \int_{\Omega \times Y_1} (|\bar{w}_1|^2 + |\bar{w}_2|^2) dx dyd\tau
  \quad \mbox{  for } t \in  [0,T].
\end{eqnarray*}
This leads to
\begin{eqnarray*}
& & \|\bar{\Lambda}_T(\bar{w}_1^1, \bar{w}_2^1) -
\bar{\Lambda}_T(\bar{w}_1^2, \bar{w}_2^2)\|_{L^2((0,T) \times \Omega \times Y_1)} \\
& \leq & \|w_1\|_{L^2((0,T) \times \Omega \times Y_1)}
+ \|w_2\|_{L^2((0,T) \times \Omega \times Y_1)} \\
& \leq &
C e^{Ct} T^{1/2} \|(\bar{w}_1^1, \bar{w}_2^1) - (\bar{w}_1^2, \bar{w}_2^2)
  \|_{L^2((0,T) \times \Omega \times Y_1)}
\end{eqnarray*}
Therefore, there exists a positive number $T_0$ such that
$\bar{\Lambda}_{T_0}$ is a contraction mapping for $0 < T_0 \leq T$. Since the choice of $T_0$ is independent of initial values, by Banach's fixed point theorem
we conclude that the problem \eqref{pd1}--\eqref{main_bc} has a solution in the sense of (S').
\hfill $\Box$

\vskip 12pt

{\it Proof of Theorem \ref{existence}. }
First, for $m > 0$
 we define  $f_{im}$, $i = 1, 2$,   $R_m$ and $Q_m$  by
\begin{eqnarray}
f_{im}(r):=
\left \{ \begin{array}{cl}
f_i(m)&\mbox{ for } r > m,\\\nonumber
f_i(r)&\mbox{ otherwise,  }
 \end{array} \right.
\quad\quad
R_m(r):=
\left \{ \begin{array}{cl}
R(m)&\mbox{ for }\;\; r > m,\\\nonumber
R(r)&\mbox{ otherwise.   }
 \end{array} \right.
\end{eqnarray}

$$
Q_m(r):=
\left \{ \begin{array}{cl}
Q(m)&\mbox{ for }\;\; r > m,\\\nonumber
Q(r)&\mbox{ for   }  \;\; |r| \leq m, \\\nonumber
Q(-m)&\mbox{ for   }  \;\; r < -m.
 \end{array} \right.
$$
Then, for each $m > 0$ by Proposition \ref{Lipschitz}  the problem \eqref{pd1}--\eqref{main_bc} with $f_1 = f_{1m}$, $f_2 = f_{2m}$,  $R = R_m$ and $Q = Q_m$ has a solution $(w_{1m}, w_{2m}, w_{3m}, w_{4m})$ in the sense of (S').

Now, for each $m$ we shall prove

(i) $ w_{1m}, w_{2m}(t) \geq 0$ a.e. on $(0,T) \times {\Omega\times Y_1}$,
$w_{3m}\geq0$ a.e. on $(0,T) \times {\Omega}$ and
 $w_{4m} \geq 0$ a.e. on $(0,T) \times \Omega\times\Gamma_1$.

In order to prove (i) we test (S2) by $w_{1m} + w_{1m}^-$, where
${\phi}^-:=  -\min\{0, \phi\}$ with $\phi^+\phi^-=0$. Then we see that
\begin{eqnarray*}
&& \frac{1}{2} \frac{d}{dt} \int_{\Omega\times Y_1}  |w_{1m}^-|^2 dxdy
+ \int_{\Omega \times Y_1} d_1 |\nabla_y w_{1m}^-|^2 dx dy  \\
& & + \int_{\Omega \times \Gamma_1} Q_m(w_4) ( \hat{R}_m(w_{1m}) - \hat{R}_m(w_{1m} + w_{1m}^-)) dx d\gamma_y
\\
 & \leq & \int_{\Omega\times Y_1}( f_{1m}(w_{1m})  - f_{2m}(w_{2m}) )w_{1m}^-dxdy \quad
 \mbox{ a.e. on } [0,T],
\end{eqnarray*}
where $\hat{R}_m$ is the primitive of ${R}_m$ with $\hat{R}_m(0) = 0$.
Note that $\hat{R}_m(w_{1m}) - \hat{R}_m(w_{1m} + w_{1m}^-) = 0 $ and
$ ( f_{1m}(w_{1m})  - f_{2m}(w_{2m}) )w_{1m}^- \leq 0$, since $f_{2m} \geq 0$ on $\mathds{R}$. Clearly,
$$ \frac{1}{2} \frac{d}{dt} \int_{\Omega\times Y_1} |w_{1m}^-|^2 dxdy
+ \int_{\Omega \times Y_1} d_1 |\nabla_y w_{1m}^-|^2 dx dy \leq 0 \quad \mbox{ a.e. on } [0,T] $$
so that $w_{1m} \geq 0$ a.e. on $(0,T) \times \Omega \times Y_1$.

Next, because $-[w_{3m}]^- \in X$, we can test (S3) by $-{w_{2m}}^-$ and  (S4) by $-w_{3m}^-$ to obtain
\begin{eqnarray}
&&  \frac{1}{2} \frac{d}{dt}
      \int_{\Omega\times Y_1}  |w_{2m}^-|^2  dx dy
  + d_2^0 \int_{\Omega \times Y_1} |\nabla_y w_{2m}^-|^2  dxdy
  + \alpha\int_{\Omega \times \Gamma_2} |w_{2m}^-|^2  dx d\gamma_y \nonumber \\
& \leq &
-  \int_{\Omega\times Y_1}( f_{1m}(w_{1m})  - f_{2m}(w_{2m}) ) w_{2m}^- dxdy
- \alpha \int_{\Omega\times\Gamma_2}  H w_{3m}{w_{2m}}^- dx d\gamma_y \nonumber \\
& \leq &
\frac{\alpha}{2}  \int_{\Omega\times\Gamma_2}  |w_{2m}^-|^2 dx d\gamma_y
 +  \frac{\alpha}{2} H^2 |\Gamma_2|  \int_{\Omega}  |w_{3m}^-|^2 dx \quad \mbox{ a.e. on } [0,T],
 \label{pro2-1}
\end{eqnarray}
and
\begin{eqnarray}
  \frac{1}{2}   \frac{d}{dt}   \int_{\Omega} |w_{3m}^-|^2 dx
  + d_3^0 \int_{\Omega} |\nabla w_{3m}^-|^2 dx
&  = & \alpha \int_{\Omega\times\Gamma_2}  (H w_{3m} - w_{2m}) w_{3m}^- dx d\gamma_y \nonumber \\
&  \leq & \alpha \int_{\Omega\times\Gamma_2}  |w_{2m}^-| |w_{3m}^-| dx d\gamma_y  \label{pro2-2}
\mbox{ a.e. on } [0,T].
\end{eqnarray}
Adding \eqref{pro2-1} and \eqref{pro2-2} and then applying Young's inequality, we get
\begin{eqnarray}
&&  \frac{1}{2} \frac{d}{dt} (
      \int_{\Omega\times Y_1}  |w_{2m}^-|^2  dx dy +  \int_{\Omega} |w_{3m}^-|^2 dx )
  + d_2^0 \int_{\Omega\times Y_1} |\nabla_y w_{2m}^-|^2  dxdy
  + d_3^0 \int_{\Omega} |\nabla w_{3m}^-|^2 dx   \nonumber \\
& \leq &
  (\frac{\alpha}{2} H^2 |\Gamma_2| + \alpha|\Gamma_2|)  \int_{\Omega}  |w_{3m}^-|^2 dx
\quad \mbox{ a.e. on } [0,T].
\nonumber
\end{eqnarray}
The application of Gronwall's inequality and the positivity of initial data give
$w_{2m} \geq 0$ a.e. on $(0,T) \times \Omega \times Y_1$ and
$w_{3m} \geq 0$ a.e. on $(0,T) \times \Omega$.

Since $\eta \geq 0$, it is easy to see that
$$ \frac{1}{2} \frac{d}{dt} \int_{\Omega \times \Gamma_1} |w_{4m}^-| dx d\gamma_y \leq 0
 \quad \mbox{ a.e. on } [0,T].
$$
Hence, we see that $w_{4m} \geq 0$ a.e. on $(0,T) \times \Omega \times \Gamma_1$.
Thus (i) is true.

Next, we shall show upper bounds of solutions as follows: To do so,
 by (A1) we can take $M_1$ and $M_2$ such that
$$  M_1 \geq  \|w_{10}\|_{L^\infty(\Omega\times Y_1)},
 M_2 \geq  \max\{ \|w_{20}\|_{L^\infty(\Omega\times Y_1)},
    H \|w_{30}\|_{L^\infty(\Omega)}, H \|w_{3}^D\|_{L^\infty(\Omega \times Y_1)},  \}
$$
and $f_1(M_1) = f_2(M_2)$. Also, we put $M_3 = \frac{M_2}{H}$,
$M_4 = \max\{ \beta_{max}, \|w_{40}\|_{L^{\infty}(\Omega \times \Gamma_1)}\}$
and $M_0 = \max\{M_1, M_2, M_3,  M_4\}$.
Then it holds:

(ii)  For any $m \geq M_0$ we have
$w_{1m}(t)\leq M_1$, $w_{2m}(t)\leq M_2$ a.e. in ${\Omega\times Y_1}$,
 $w_{3m}(t)\leq M_3$ a.e. in ${\Omega}$ and
$w_{4m}(t)\leq M_4$ a.e. on $\Omega\times\Gamma_1$ for a.e. $t \in [0,T]$.

In fact, let $m \geq M_0$ and consider $w_{1m} - (w_{1m}-M_1)^+$, $(w_{2m}-M_2)^+$ and $(w_{3m}-M_3)^+$ as  test functions in (S2) $\sim$ (S4). Then we observe that
\begin{eqnarray}
&& \frac{1}{2} \frac{d}{dt}
 \int_{\Omega\times Y_1}  |(w_{1m}- M_1)^+|^2 dx dy
+ d_1^0  \int_{\Omega\times Y_1} |\nabla_y(w_{1m}-M_1)^+|^2 dxdy \nonumber  \\
& & +  \int_{\Omega\times\Gamma_1} Q_m(w_{4m})
  (\hat{R}_m(w_{1m}) - \hat{R}_m(w_{1m} - (w_{1m}-M_1)^+))
 dx d\gamma_y  \nonumber
 \\
&\leq &
 \int_{\Omega\times Y_1} (- f_{1m}(w_{1m}) + f_{2m}(w_{2m}) ) (w_{1m}- M_1)^+ dxdy, \label{11}
\end{eqnarray}
\begin{eqnarray}
 &&  \frac{1}{2} \frac{d}{dt}  \int_{\Omega\times Y_1}  |(w_{2m}-M_2)^+|^2 dx dy
+d_2^0   \int_{\Omega\times Y_1} |\nabla_y(w_{2m}-M_2)^+|^2 dxdy \nonumber \\
&\leq & \!\!\!
\int_{\Omega\times Y_1}\!\!\!(f_{1m}(w_{1m})   - f_{2m}(w_{2m}))(w_{2m}-M_2)^+dxdy   \label{12} \\
& &  + \! \alpha \!\! \int_{\Omega\times\Gamma_2}\!\!\!\!\! (H w_{3m} - w_{2m}) (w_{2m}-M_2)^+ dx
 d\gamma_y,  \nonumber
\end{eqnarray}
\begin{eqnarray}
&& \frac{1}{2} \frac{d}{dt} \int_{\Omega}
           |(w_{3m}-M_3)^+|^2 dx  +
d_3^0  \int_{\Omega} |\nabla(w_{3m} - M_3)^+|^2 ) dx \nonumber \\
&  \leq &
-\alpha \int_{\Omega \Gamma_2} (Hw_{3m} - w_{2m}) (w_{3m}-M_3)^+ dxd\gamma_y
  \quad \mbox{ a.e. on } [0,T]. \label{14}
\end{eqnarray}
Here, we note that $\hat{R}_m(w_{1m}) - \hat{R}_m(w_{1m} - (w_{1m}-M_1)^+) \geq 0$. Adding \eqref{11}--\eqref{14}, we get
\begin{eqnarray}
&& \frac{1}{2} \frac{d}{dt} \left(\int_{\Omega\times Y_1}
  (|(w_{1m}-M_1)^+|^2 +  |(w_{2m}-M_2)^+|^2)dxdy
 + \int_{\Omega} |(w_{3m}-M_3)^+|^2  dx \right) \nonumber \\
&& +  \int_{\Omega\times Y_1} \left(d_1^0 |\nabla(w_{1m} - M_1)^+|^2
 + d_2^0 |\nabla(w_{2m}-M_2)^+|^2  \right) dxdy
+ d_3^0 \int_{\Omega}  |\nabla(w_{3m}-M_3)^+|^2 dx \nonumber \\
 &\leq &
 \int_{\Omega\times Y_1}
  (- f_{1m}(w_{1m}) +  f_{2m}(w_{2m})) ((w_{1m}-M_1)^+ - (w_{2m}-M_2)^+) dxdy \label{ine}
 \\&&+
 \alpha \int_{\Omega\times\Gamma_2} ((Hw_{3m}-w_{2m})(w_{2m}-M_2)^+
 +(w_{2m}-{H}w_{3m})(w_{3m}-M_3)^+)dxd\gamma_y  \mbox{ a.e. on } [0,T]. \nonumber
\end{eqnarray}
We estimate the first term on the r.h.s of \eqref{ine} by making use of $f_{1m}(M_1)=f_{2m}(M_2)$
and  the Lipschitz continuity of $f_{im}$, $i = 1, 2$, as follows: We have
\begin{eqnarray*}
&&
 \int_{\Omega\times Y_1} (- f_{1m}(w_{1m}) + f_{1m}(M_{1}) - f_{2m}(M_{2})
+ f_{2m}(w_{2m})) (w_{1m}-M_{1})^+dxdy
\\
&&
 + \int_{\Omega\times Y_1}( f_{m1}(w_{1m})-f_{1m}(M_{1})+f_{2m}(M_{2}) - f_{2m}(w_{2}))
    (w_{2m}-M_2)^+dxdy
\\
& \leq &
 \int_{\Omega\times Y_1} (f_{2m}(w_{2m})-f_{2m}(M_{2})) (w_{1m}-M_{1})^+ dxdy \\
& & + \int_{\Omega\times Y_1} (f_{1m}(w_{1m})-f_{1m}(M_{1})) (w_{2m}-M_{2})^+dxdy \\
& \leq &
 C \int_{\Omega\times Y_1}  (|(w_{2m}-M_2)^+|^2 +|(w_{1m}-M_1)^+|^2)dxdy
 \quad \mbox{ a.e. on } [0,T].
\end{eqnarray*}

We estimate the second term on the r.h.s in \eqref{ine} as follows:
\begin{eqnarray*}
&& \alpha  \int_{\Omega\times\Gamma_2}
       (Hw_{3m}-HM_3+M_2-w_{2m})(w_{2m}-M_2)^+ dxd\gamma_y \\
&& +\alpha \int_{\Omega\times\Gamma_2} (w_{2m}-M_2+{H}{}(M_3-w_{3m}))(w_{3m}-M_3)^+
           dxd\gamma_y
\\
&\leq&
\alpha H \int_{\Omega\times\Gamma_2} (w_{m3} -  M_3)(w_{2m} -  M_2)^+ dx d\gamma_y
- \alpha  \int_{\Omega\times\Gamma_2} |(w_{2m} -  M_2)^+|^2 dx d\gamma_y
\\ &&
+ \alpha \int_{\Omega\times\Gamma_2} (w_{2m} -  M_2)(w_{3m} -  M_3)^+ dx d\gamma_y
- \alpha H \int_{\Omega\times\Gamma_2} |(w_{3m} -  M_3)^+|^2  dx d\gamma_y
\\
&\leq&
(\alpha H^2 + \alpha) \int_{\Omega\times\Gamma_2} |(w_{3m} -  M_3)^+|^2  dx d\gamma_y
\quad \mbox{ a.e. on } [0,T].
\end{eqnarray*}
Now, \eqref{ine} becomes
\begin{eqnarray*}
&& \frac{1}{2} \frac{d}{dt} \left(\int_{\Omega\times Y_1}
  (|(w_{1m}-M_1)^+|^2 +  |(w_{2m}-M_2)^+|^2)dxdy
 + \int_{\Omega} |(w_{3m}-M_3)^+|^2  dx \right) \nonumber \\
&& +  \int_{\Omega\times Y_1} \left(d_1^0 |\nabla(w_{1m} - M_1)^+|^2
 + d_2^0 |\nabla(w_{2m}-M_2)^+|^2  \right) dxdy \nonumber \\
& & + d_3^0 \int_{\Omega}  |\nabla(w_{3m}-M_3)^+|^2 dx \nonumber \\
&\leq&
 C \int_{\Omega\times Y_1} (|(w_{2m}-M_2)^+|^2+| (w_{1m}-M_1)^+|^2) dxdy \\
& & +C \int_{\Omega}|(w_{3m}-M_3)^+|^2 dx  \quad \mbox{ a.e. on } [0,T].
\end{eqnarray*}
Applying Gronwall's inequality, we get
$$
 \int_{\Omega\times Y_1} (|(w_{1m}(t)-M_1)^+|^2  +
|(w_{2m}(t)-M_2)^+|^2)dxdy +
\int\limits_{\Omega} |(w_{3m}(t)-M_3)^+|^2  dx \leq 0
\mbox{ for } t \geq 0.
$$
Hence, $w_{1m}\leq M_1, w_{2m}\leq M_2$ a.e. in $\Omega\times Y_1$ and $ w_{3m}\leq M_3$ a.e. in $\Omega$ for $t\in(0,T)$.

To show that $w_{4m}$ is bounded on $\Omega\times\Gamma_1$, we test \eqref{od1} with
$(w_{4m}- M_4)^+$ and using (A2) leads to
\begin{eqnarray*}\label{}
& &  \frac{1}{2} \frac{d}{dt} \int\limits_{\Omega\times \Gamma_1}
   |(w_{4m}- M_4)^+|^2 dxd\gamma_y \\
&  = &
 \int\limits_{\Omega\times \Gamma_1} R_m(w_1)Q_m(w_4) (w_{4m}- M_4)^+ dxd\gamma_y
 \leq  0  \quad \mbox{ a.e. on } [0,T].
\end{eqnarray*}
This shows that $w_{4m}\leq M_4$ a.e on $(0,T) \times \Omega\times\Gamma_1$.
Thus we have (ii).

Accordingly,  by (i) and (ii)   $(w_{1m}, w_{2m}, w_{3m}, w_{4m})$ satisfies
the conditions (S1) $\sim$ (S4) for $m \geq M_0$.
Thus we have proved this theorem. \hfill $\Box$

\vskip 12pt
{\it Proof of Theorem \ref{Uniqueness}. } Let $(w_{1j}, w_{2j}, w_{3j}, w_{4j})$, $j = 1, 2$, be solutions \eqref{pd1}--\eqref{main_bc} satisfying (S1) $\sim$ (S4).
Since all $w_{ij}$, $i = 1, 2, 3, 4$, $j = 1, 2$,   are bounded,  $(w_{1j}, w_{2j}, w_{3j}, w_{4j})$ is also a solution of \eqref{pd1}--\eqref{main_bc} with $f_1 = f_{1m}$, $f_2 = f_{2m}$, $R = R_m$ and $Q = Q_m$ for some positive constant $m$.  Then Proposition \ref{Lipschitz} guarantees the uniqueness.
This proves the conclusion of Theorem  \ref{Uniqueness}. \hfill $\Box$

\section*{Acknowledgements}
 We acknowledge fruitful discussions on this subject with M. Neuss-Radu (Erlangen), O. Lakkis (Sussex), and V. Chalupecky (Fukuoka).
  A. M. and T. A. thank both science foundations NWO and JSPS for supporting financially the
   Dutch-Japanese seminar 'Analysis of non-equilibrium evolution problems, selected
   topics in material and life sciences', during which this paper was completed.

\bibliographystyle{plain}

\end{document}